# Level crossings and other level functionals of stationary Gaussian processes*

## Marie F. Kratz

*ESSEC*
*MAP5, UMR 8145, Université Paris V*
*SAMOS-MATISSE, UMR 8174, Université Paris I*
*e-mail:* `kratz@essec.fr`

**Abstract:** This paper presents a synthesis on the mathematical work done on level crossings of stationary Gaussian processes, with some extensions. The main results [(factorial) moments, representation into the Wiener Chaos, asymptotic results, rate of convergence, local time and number of crossings] are described, as well as the different approaches [normal comparison method, Rice method, Stein-Chen method, a general $m$-dependent method] used to obtain them; these methods are also very useful in the general context of Gaussian fields. Finally some extensions [time occupation functionals, number of maxima in an interval, process indexed by a bidimensional set] are proposed, illustrating the generality of the methods. A large inventory of papers and books on the subject ends the survey.

## Contents

---

*This is an original survey paper

## 1. Introduction

This review presents the mathematical aspect of the work done on level crossings and upcrossings of a level, and insists on the different approaches used to tackle a given problem. We speak only briefly about examples or applications which of course were the reason for such theoretical work and which have been greatly developing in the extreme value field over the past years. This review of mathematical results and methods can also be very useful in mathematical statistics applications, reliability theory and other areas of applications of stochastic processes.

When possible, we choose to privilege conditions on the behavior of spectral or covariance functions rather than general conditions of mixing, since they are easier to work with.

We mainly consider continuous parameter processes, since they appear in most of the mathematical models that describe physical phenomena. Nevertheless one approach to work on level crossings or on upcrossings of a level that are analogous to the exceedances used in the discrete case, can be through discretization, using discrete parameter results.

Let $X = (X_t, t \in \mathbb{R}^d)$ be a real stochastic process. Our main concern is the measure of the random set of level $x$ defined by $C_x^X := \{t : X_t = x\}$, by taking the same notation as Wschebor (see [163]). The unidimensional case will be exposed in detail with results and methods (especially the ones which can/could be easily adapted to a higher dimension), whereas only a partial view of the multidimensional case will be given.

As a preliminary, let us mention two important results (see [31] for the first one and for instance [90] for the second one) which reinforce our interest in the study of $C_x^X$.

**Theorem 1.1** *(Bulinskaya, 1961)*
*Let $X = (X_s, 0 \le s \le t)$ be an a.s. continuously differentiable stochastic process with one-dimensional density $p_s(u)$ bounded in $u$ for all $s$.*
*Then for any level $x$ the probability of existence of a point $s$ such that the events $(X(s) = x)$ and $(\dot{X}(s) = 0)$ occur simultaneously, is equal to 0. In particular, the probability that $X_t$ becomes tangent to the level $x$, is equal to 0, in other*

*words the probability of contingence of the level x by the process X is equal to 0.*

In the case of a Gaussian process, as a consequence of the separability property and the Tsyrelson theorem, we have

**Theorem 1.2** *For any Gaussian process X on a arbitrary parameter set T, for any x such that*

$$I\!\!P[\sup_{s \in T} X(s) \ < \ x] > 0,$$

*there is no contingence, with probability one.*

Let us give the notations and specific symbols used through this survey.

$L^p(\Omega) := \{X : I\!\!E[|X|^p] < \infty\}; \ ||.||_p : L^p$-norm.

$\mathbb{1}_A$: indicator function of a set $A$.

$\phi_{s,t}(x, \dot{x}, y, \dot{y})$: density function of $(X_s, \dot{X}_s, X_t, \dot{X}_t)$;

$\phi_t(x, \dot{x}, y, \dot{y}) := \phi_{0,t}(x, \dot{x}, y, \dot{y})$.

r.v.: random variable

CLT: central limit theorem

s.t.: such that

$\xrightarrow{d}$ stands for convergence in distribution.

## 2. Crossings of Gaussian processes

Studies on level-crossings by stationary Gaussian processes began about sixty years ago. Different approaches have been proposed. Here is a survey of the literature on the number of crossings of a given level or of a differentiable curve in a fixed time interval by a continuous spectrum Gaussian process.

Besides the well-known books about Extremes, let us quote a short survey by Slud in 1994 (see [148] or [149]), a more general survey about extremes including a short section about level crossings by Leadbetter et al. in 1988 (see [86]) and by Rootzen in 1995 (see [132]), and another one by Piterbarg in his 1996 book (see [122], in particular for some methods described in more detail than here). Our purpose here is to focus only on crossing counts in order to be more explicit about the subject, not only recalling the main results, but also giving the main ideas about the different methods used to establish them. To make the methods easier to understand, we consider mainly the crossing counts of a given level. Note that the problem of curve ($\psi$) crossings by a stationary process $X$ may also be regarded as a zero-crossing problem for the non stationary process $X^* := X - \psi$ (but stationary in the sense of the covariance), as pointed out by Cramér and Leadbetter (see [38]).

Let $X = (X_s, s \geq 0)$ be a real stationary Gaussian process with variance one and a.s. continuous sample functions.

Let $\psi$ be a continuously differentiable function.

We denote $N_I(x)$ or $N_I(\psi)$ the number of crossings of a given level $x$ or of a curve $\psi(.)$ respectively, $X$ on the interval $I$, and $N_I^+(x)$ the number of up-crossings of $x$ by $X$ (recall that $X_s$ is said to have an upcrossing of $x$ at $s_0 > 0$ if for some $\varepsilon > 0$, $X_s \le x$ in $(s_0 - \varepsilon, s_0)$ and $X_s \ge x$ in $(s_0, s_0 + \varepsilon)$). Hereafter this will be simply denoted $N_t(x)$ or $N_t(\psi)$ or $N_t^+(x)$ respectively, when $I = [0, t]$, $t > 0$.

## 2.1. Moments and factorial moments

Distributional results about level or curve crossings by a Gaussian process are often obtained in terms of factorial moments for the number of crossings, that is why many authors have worked not only on moments but also on factorial moments to find out expressions or conditions (in terms of the covariance function of the process) for their finiteness, since a certain number of applications require to know if they are finite.

Note also that the conditions governing finite crossing moments are local ones, since Hölder inequality implies that $(\mathbb{E}[N_{2t}(x)])^k \le 2^k(\mathbb{E}[N_t(x)])^k$, and therefore when $(\mathbb{E}[N_t(x)])^k$ is finite for some $t$, it is finite for all $t$.

### 2.1.1. Introduction

● *Kac's method and formula (1943).*
We can start with Kac (see [73]), who studied the number $N$ of zeros of a Gaussian random polynomial on some bounded interval of $\mathbb{R}$. To compute $\mathbb{E}[N]$, he proposed a method which used, so to speak, in a formal way, the approximation of the "Dirac function" $\delta_0(x)$ by the function $\dfrac{1}{2\varepsilon}\mathbb{1}_{[-\varepsilon, \varepsilon)}(x)$. He gave the first heuristic expression of $N_t(0)$ as a function of the process $X$:

**Theorem 2.1** *(Kac formula, 1943)*

$$N_t(0) = \frac{1}{2\pi}\int_{-\infty}^{\infty}\int_0^t cos(\zeta X_s)|\dot{X}_s|ds d\zeta, \quad \text{with probability one.}$$

**From now on, we suppose w.l.o.g. that X is centered**, with correlation function $r(t) = \mathbb{E}[X_0 X_t]$, given also in respect with the spectral function $F$ as $r(t) = \displaystyle\int_0^{\infty}\cos(\lambda t)dF(\lambda)$. We denote $\lambda_2$ the second moment (when it exists) of the spectral function, i.e. $\lambda_2 = \displaystyle\int_0^{\infty}\lambda^2 dF(\lambda)$.

● *Rice (1945).*
One of the best known first results is the one of Rice (see [129]) who proved by intuitive methods related to those used later by Cramér and Leadbetter (as the discretization method described below), that

**Theorem 2.2** *(Rice formula, 1945)*

$$\mathbb{E}[N_t(x)] = te^{-x^2/2} \sqrt{-r''(0)} / \pi. \tag{2.1}$$

It means that the mean number of crossings is the most important for zero crossings, and decreases exponentially with the level.

• *Itô (1964), Ylvisaker (1965).*
Itô and Ylvisaker have proved (2.1) and therefore that the mean number of crossings is finite exactly when the second spectral moment is finite:

**Theorem 2.3** *(Itô, 1964; Ylvisaker, 1965)*

$$\mathbb{E}[N_t(x)] < \infty \quad \Leftrightarrow \quad \lambda_2 < \infty \quad \Leftrightarrow \quad -r''(0) < \infty.$$

• *The discretization method.*
From the intuitive method developed by Rice at the beginning of the 40s (see [129]), Ivanov in 1960 (see [70]), Bulinskaya in 1961 (see [31]), Itô in 1964 (see [71]) and Ylvisaker in 1965 (see [165]) proposed rigorous proofs for zero countings. There followed the general formulation due to Leadbetter in 1966 (see [38] p.195 or [85] p.148), known as the method of discretization, which also applies to non-normal processes.
This method is based on the approximation of the continuous process $X = (X_t, t \geq 0)$ by the sequence $(X(jq), j = 1, 2, \cdots)$ where $q$ satisfies some conditions related to the level $x$ of upcrossings:

$$q = q(x) \to 0 \quad \text{and} \quad xq \downarrow 0 \quad \text{as} \quad x = x(t) \to \infty \quad \text{(or equivalently as } t \to \infty). \tag{2.2}$$

Consider any sequence $(q_n)$ s.t. $q_n \downarrow 0$ and the number $N_{n,I}$ of points $(jq_n, j = 1, 2, \cdots)$, in a fixed and bounded interval I, s.t. $\begin{cases} (j-1)q_n \in I \text{ and } jq_n \in I, \\ X((j-1)q_n) < x < X(jq_n). \end{cases}$
It can be proved (see [38] p.195, and [85], lemma 7.7.2) that

$$N_{n,I} \to N_I^+(x) \ a.s. \ as \ n \to \infty, \tag{2.3}$$

hence that $N_I^+(x)$ is a (possibly infinite-valued) r.v., and that (see [85])
$\mathbb{E}[N_{n,I}] \underset{n \to \infty}{\to} \mathbb{E}[N_I^+(x)]$.
Following the notation of Leadbetter et al. (see [85]), let

$$J_q(x) = \frac{1}{q} \mathbb{P}[X_0 < x < X_q], \quad q > 0. \tag{2.4}$$

By choosing $I = (0, t]$, we have $\nu_n := t/q_n$ points $jq_n \in I$ and by stationarity, it yields $\mathbb{E}[N_n] = (\nu_n - 1)\mathbb{P}[X_0 < x < X_{q_n}] \sim tJ_{q_n}(x)$; since the sequence $(q_n)$ is arbitrary, it implies

$$\mathbb{E}[N_t^+(x)] = t \lim_{q \downarrow 0} J_q(x) \tag{2.5}$$

and hence that

$$\mathbb{E}[N_t^+(x)] = t\,\mathbb{E}[N_1^+(x)].\qquad(2.6)$$

Under mild conditions, $\lim_{q\to 0} J_q(x)$, and therefore $\mathbb{E}[N_1^+(x)]$, can be expressed in a simple integral form (by writing $J_q$ as $\mathbb{P}[X_0 > x \ , \ \tilde{X}_q > q^{-1}(x - X_0)]$, with $\tilde{X}_q := q^{-1}(X_q - X_0)$, and by doing some change of variable), giving in the normal case (see [85], lemma 7.3.1):

**Lemma 2.1** *(Leadbetter et al.)*
*Let $X = (X_t, t \geq 0)$ be a standardized stationary normal process with its second spectral moment $\lambda_2$ finite and let $(q = q(x), x)$ satisfy (2.2). Then, as $q \to 0$,*

$$J_q(x) \ \sim \ \int_0^\infty z p(x, z) dz = \frac{\lambda_2^{1/2}}{2\pi} e^{-x^2/2},$$

*$p(.,.)$ denoting the joint density of $(X, \dot{X})$ and $J_q$ satisfying (2.4).*

Therefore this last result, together with (2.5), provide the Rice formula for the number of upcrossings of any fixed level $x$ per unit time by a standardized stationary normal process, namely

$$\mathbb{E}[N_1^+(x)] = e^{-x^2/2} \ \lambda_2^{1/2}/(2\pi).$$

• *Cramér and Leadbetter (1965,1967), Ylvisaker (1966).*
In the 60s, following the work of Cramér, generalizations to curve crossings and higher order moments for $N_t(.)$ were considered in a series of papers by Cramér and Leadbetter (see [38]) and Ylvisaker (see [166]).
As regards curve crossings, the generalized Rice formula was obtained:

**Lemma 2.2** *(Generalized Rice formula; Ylvisaker, 1966; Cramér and Leadbetter, 1967).*
*Let $X = (X_t, t \geq 0)$ be a mean zero, variance one, stationary Gaussian process, with twice differentiable covariance function $r$. Let $\psi$ be a continuous differentiable real function. Then*

$$\mathbb{E}[N_t(\psi)] = \sqrt{-r''(0)} \int_o^t \varphi(\psi(y)) \left[ 2\varphi\left(\frac{\psi'(y)}{\sqrt{-r''(0)}}\right) + \frac{\psi'(y)}{\sqrt{-r''(0)}}\left(2\Phi\left(\frac{\psi'(y)}{\sqrt{-r''(0)}}\right) - 1\right)\right] dy,\ (2.7)$$

*where $\varphi$ and $\Phi$ are the standard normal density and distribution function respectively.*

Again the authors used the discretization method and approximated the continuous time number $N_t(\psi)$ of $\psi$-crossings by the discrete-time numbers $(N_\psi(t, q_n))$ of crossings of continuous polygonal curves agreeing with $\psi$ at points $jq_n$ $(j = 0, 1, \cdots, 2^n)$, with time steps $q_n = t/2^n$:

$$N_\psi(t, q_n) := \sum_{j=0}^{2^n-1} \mathbb{1}_{\left((X_{jq_n} - \psi(jq_n))(X_{(j+1)q_n} - \psi((j+1)q_n)) < 0\right)}\qquad(2.8)$$

to obtain that $N_\psi(t, q_n) \uparrow N_t(\psi)$ with probability one, as $n \to \infty$. $\square$

### 2.1.2. Moments and factorial moments of order 2

$\bullet$ Cramér and Leadbetter proposed a sufficient condition (known nowdays as the *Geman condition*) on the correlation function of $X$ (stationary case) in order to have the random variable $N_t(x)$ belonging to $L^2(\Omega)$, namely

**Theorem 2.4** *(Cramér et al., 1967)*

$$If \quad \exists \delta > 0, \quad L(t) := \frac{r''(t) - r''(0)}{t} \in L^1([0, \delta], dx) \qquad (2.9)$$

$$then \quad \mathbb{E}[N_t^2(0)] < \infty.$$

$\bullet$ Explicit (factorial) moment formulas for the number of crossings have been obtained by Cramér and Leadbetter (see [38]), Belyaev (see [19]) and Ylvisaker (see [165]), based on careful computations involving joint densities of values and derivatives of the underlying Gaussian process.

In particular, the second moment of $N_t(x)$ is given by

$$\mathbb{E}[N_t^2(x)] = \mathbb{E}[N_t(x)] + 2 \int_0^t (t - u) \int_{\mathbb{R}^2} |\dot{x}||\dot{y}| \phi_u(x, \dot{x}, x, \dot{y}) d\dot{x} d\dot{y} du.$$

Concerning the second factorial moment, Cramér & Leadbetter provided an explicit formula for the number of zeros of the process $X$ (see [38], p.209), from which the following formula for the second factorial moment of the number of crossings of a continuous differentiable real function $\psi$ by $X$ can be deduced:

$$\mathbb{E}[N_t(\psi)(N_t(\psi) - 1)] = \int_0^t \int_0^t \int_{\mathbb{R}^2} |\dot{x}_1 - \dot{\psi}_{t_1}||\dot{x}_2 - \dot{\psi}_{t_2}| \phi_{t_1, t_2}(\psi_{t_1}, \dot{x}_1, \psi_{t_2}, \dot{x}_2) d\dot{x}_1 d\dot{x}_2 dt_1 dt_2 \quad (2.10)$$

where the density $\phi_{t_1, t_2}$ is supposed non-singular for all $t_1 \neq t_2$. The formula holds whether the second factorial moment is finite or not.

$\bullet$ *Ershov (1967).*
This author proved (see [51]) that whenever a Gaussian stationary process $X = (X_t, t \in \mathbb{R})$ with covariance function $r$ is with mixing (i.e. that $|r(|i - j|)| \leq f(|i - j|)$ with $\lim_{k \to \infty} f(k) = 0$), then the number of its $x$-upcrossings on all $\mathbb{R}$ cannot hold a finite second moment:

**Theorem 2.5** *(Ershov, 1967)*
If $r(t) \to 0$ as $t \to \infty$,
then $var\left(N_t^+(x)\right) \to \infty$, as $t \to \infty$.

$\bullet$ *Geman (1972).*
Geman proved in 1972 (see [54]) that the condition (2.9) was not only sufficient but also necessary, by showing that if $\frac{r''(t) - r''(0)}{t}$ diverges on $(0, \delta)$ then so

does the integral
$\int_{-\infty}^{\infty} \int_{-\infty}^{\infty} |xy| \phi_s(0, x, 0, y) dx dy$ appearing in the computation of the second moment. Thus we have

**Theorem 2.6** *(Geman condition, 1972)*

$$\text{Geman condition } (2.9) \quad \Leftrightarrow \quad I\!\!E[N_t^2(0)] < \infty.$$

Of course, this result holds for any stationary Gaussian process when considering the number of crossings of the mean of the process.

- *Wschebor (1985).*
This author provided in [163] (under the Geman condition) an explicit expression in the case of two different levels $x$ and $y$, with $x \neq y$, namely:

$$I\!\!E[N_t(x)N_t(y)] = \int_0^t (t-u) \int_{I\!\!R^2} |\dot{x}||\dot{y}| \left( \phi_u(x, \dot{x}, y, \dot{y}) + \phi_u(y, \dot{x}, x, \dot{y}) \right) d\dot{x} d\dot{y} du.$$

Note that this expression differs from the one of Cramér et al. when $x \to y$, which means that the function $I\!\!E[N_t(x)N_t(y)]$ is not continuous on the diagonal.

- *Piterbarg (1982), Kratz and Rootzén (1997).*
Let us assume that the correlation function $r$ of $X$ satisfies

$$r(s) = 1 - \frac{s^2}{2} + o(s^2) \qquad \text{as} \quad s \to 0, \tag{2.11}$$

$$I\!\!E\{(X'(s) - X'(0))^2\} = 2(r''(s) - r''(0)) \leq c|s|^\gamma, \qquad s \geq 0, \quad 0 < \gamma \leq 2, \tag{2.12}$$

that $r(s)$ and its first derivative decay polynomially,

$$|r(s)| \leq C s^{-\alpha}, \qquad |r(s)| + r'(s)^2 \leq C_0 s^{-\alpha}, \qquad s \geq 0, \tag{2.13}$$

for some $\alpha > 2$ and constants $c, C, C_0$,
and that the range of $t, x$ is such that

$$\varepsilon \leq I\!\!E[N_x^+(t)] \leq K_0 \tag{2.14}$$

for some fixed $K_0, \varepsilon > 0$; by using Piterbarg's notations, let

$$m = m(s) = \frac{r'(s)}{1 + r(s)}, \quad \sigma = \sigma(s) = \sqrt{\frac{1 - r(s)^2 - r'(s)^2}{1 - r(s)^2}}$$

and

$$m^*(s) = \sup_{v \geq s > 0} \frac{m^2(v)}{2\pi \sqrt{1 - r^2(v)}}, \quad \sigma^* = \sup_{s \geq 0} \frac{\sigma^2(s)}{2\pi \sqrt{1 - r^2(s)}}.$$

Piterbarg (see [121]) proposed some bounds of the second factorial moment of the number of upcrossings; Kratz and Rootzén (see [84]) gave a small variation of his result, but with more precise bounds, namely

**Lemma 2.3** *(Kratz and Rootzén, 1997)*
*Suppose that the previous hypotheses are satisfied (with $\gamma = 2$).*
*Then, for $t, x \geq 1$,*

$$\mathbb{E}\left[N_x^+(t)(N_x^+(t) - 1)\right] \leq K_3 t x^{2+2/\alpha} e^{-\frac{x^2}{2}\left(1 + \inf_{s \geq 0} \rho(s)\right)} + K_4 t^2 x^2 e^{-x^2}, \quad (2.15)$$

*with $\rho(.)$ defined in (2.42) below, $K_3 = (3.1\sigma^* + 2.3m^*(0))(3C_0)^{1/\alpha}$, $K_4 = 2.6\sigma^* + 2.3m^*(0)$.*

• *Kratz and León (2006).*
Recently, Kratz and León (see [83]) have succeded in generalizing theorem 2.6 to any given level $x$ (not only the mean of the process) and to some differentiable curve $\psi$.
Note that the problem of finding a simple necessary and sufficient condition for the number of crossings of any level has already been broached in some very interesting papers that proposed sufficient conditions as the ones of Cuzick (see [40], [43], [44]). But getting necessary conditions remained an open problem for many years. At present, the solution of this problem is enunciated in the following theorem.

**Theorem 2.7** *(Kratz and León, 2006)*

*1) For any given level $x$, we have*

$$\mathbb{E}[N_t^2(x)] < \infty \iff \exists \delta > 0, \ L(t) = \frac{r''(t) - r''(0)}{t} \in L^1([0, \delta], dx) \ \textit{(Geman condition)}.$$

*2) Suppose that the continuous differentiable real function $\psi$ satisfies for some $\delta > 0$,*

$$\int_0^\delta \frac{\gamma(s)}{s} ds < \infty, \quad \textit{where } \gamma(\tau) \textit{ is the modulus of continuity of } \dot{\psi}.$$

*Then*

$$\mathbb{E}[N_t^2(\psi)] < \infty \iff L(t) \in L^1([0, \delta], dx).$$

This smooth condition on $\psi$ is satisfied by a large class of functions which includes in particular functions whose derivatives are Hölder.

The method used to prove that the Geman condition keeps being the sufficient and necessary condition to have a finite second moment is quite simple.
It relies mainly on the study of some functions of $r$ and its derivatives at the neighborhood of 0, and the chaos expansion of the second moment, a notion which will be explained later.

*2.1.3. Factorial moments and moments of higher order*

• Concerning moments of order higher than 2, Cramér and Leadbetter (see [38]) got in the stationary case, under very mild conditions, results that Belyaev (see [19]) derived in the non-stationary case under slightly more restrictive conditions, weakened by Ylvisaker (see [166]) in the stationary case but which may also be adapted to cover nonstationary cases. Let us give for instance the $k$th factorial moment of the number of zero crossings by a stationary Gaussian process:

$$M_t^k(0) := \; \mathbb{E}[N_t(0)(N_t(0) - 1) \cdots (N_t(0) - k + 1)]$$
$$= \int_0^t dt_1 \cdots \int_0^t dt_k \mathbb{E}\left[ \prod_{i=1}^k |\dot{X}(t_j)| \mid X(t_j) = 0, \; 0 \le j \le k \right] p(0, ..., 0) \; (2.16)$$

where $p(x_1, \cdots, x_k)$ is the joint density of the r.v. $(X_{t_1}, \cdots, X_{t_k})$.

• *Belyaev (1967).*
Concurrently to this study, Belyaev proposed in [19] a sufficient condition for the finiteness of the $k$th factorial moment $M_t^k(0)$ for the number $N_t(0)$ of zero crossings on the interval $[0, t]$ in terms of the covariance matrix $\Sigma_k$ of $(X_{t_1}, \cdots, X_{t_k})$ and of

$$\sigma_i^2 := var\left( \dot{X}_{t_i} \mid X_{t_j} = 0, \; 1 \le j \le k \right); \tag{2.17}$$

**Theorem 2.8** *(Belyaev, 1967)*

$$If \quad \int_0^t dt_1 \cdots \int_0^t dt_k \left( \frac{\prod_{i=1}^k \sigma_i^2}{det\Sigma_k} \right)^{1/2} < \infty \tag{2.18}$$

$$then \quad M_t^k(0) := \mathbb{E}[N_t(0)(N_t(0) - 1) \cdots (N_t(0) - k + 1)] < \infty.$$

• *Cuzick (1975-1978).*
Cuzick proved in 1975 (see [40]) that Belyaev's condition (2.18) for the finiteness of the $k$th factorial moments for the number of zero crossings, was not only sufficient but also necessary:

**Theorem 2.9** *(Cuzick condition, 1975)*

$$Condition \; (2.18) \quad \Leftrightarrow \quad \int_0^\varepsilon d\Delta_1 \cdots \int_0^\varepsilon d\Delta_{k-1} \left( \frac{\prod_{i=1}^k \sigma_i^2}{det\Sigma_k} \right)^{1/2} < \infty \tag{2.19}$$
$$for \; some \; \varepsilon > 0, \quad where \; \sigma_i \; satisfies \; (2.17)$$
$$\Leftrightarrow \quad M_t^k(0) = \mathbb{E}[N_t(0)(N_t(0) - 1) \cdots (N_t(0) - k + 1)] < \infty.$$

The proof that condition (2.18) is equivalent to condition (2.19) is immediate with the change of variables $\Delta_i := t_{i+1} - t_i$, $1 \le i \le k - 1$, after noticing the symmetry of the integrand in (2.18).

Now the necessity of (2.19) comes from the fact that the lemma given below implies that $(2.16) > C \int_0^t dt_1 \cdots \int_0^t dt_k \left( \frac{\prod_{i=1}^k \sigma_i^2}{det \Sigma_k} \right)^{1/2}$, with $C$ some constant.

Then Cuzick tried to derive from his result (2.19) simpler sufficient conditions to have the finiteness of the $k$th (factorial) moments for the number of crossings. In particular, in his 75s paper (see [40]), he proved that $M_t^k(0) < \infty$ for all $k$, for a covariance function $r$ having a behavior near 0 such that

$r(s) = 1 - \frac{s^2}{2} + c\frac{s^3}{3!} + o(s^3)$ as $s \to 0$ (with $c$ some constant $> 0$).

Later (see [43] and [44]), for $X$ with path continuous $n$th derivative $X^{(n)}$ and spectral distribution function $F(\lambda)$, he proposed a series of sufficient conditions involving $F$ and $\sigma_n^2(h) := I\!\!E[(X_{t+h}^{(n)} - X_t^{(n)})^2]$ for having $I\!\!E[(N_t(0))^k] < \infty$.
Let us give an example in terms of the spectral density of $X$ among the sufficient conditions he proposed.

**Theorem 2.10** *(Cuzick, 1978)*
*If $X$ has a spectral density given by $f(\lambda) := 1/\left(1 + \lambda^3 |\log \lambda|^\alpha\right)$, then for $k \geq 2$ and $\alpha > 3k/2 - 1$, $M_t^k(0) < \infty$.*

Those results are not restricted to zero crossings and would also apply to a large family of curves (see [39] and [43]).
However, necessary conditions are more difficult for higher moments, the main difficulty lying in obtaining sharp lower bounds for the $\sigma_i^2$ defined in (2.19).

● *Marcus (1977): generalized Rice functions.*
Marcus (see [101]) generalized results of Cramér and Leadbetter and Ylvisaker by considering not only Gaussian processes but also by computing quantities such as $I\!\!E[N_t^{j_1}(x_1) \cdots N_t^{j_k}(x_k)]$ for levels $x_1, \cdots, x_k$ and integers $j_1, \cdots, j_k$, called generalized Rice functions.
For the proofs, the author returns to the approach used by Kac (see [73]), Ivanov (see [70]) and Itô (see [71]) to obtain the mean number of crossings at a fixed level, which consists in first finding a function that counts the level crossings of a real valued function, then in substituting $X$ for the function and finally, in considering the expectation.

● *Nualart and Wschebor (1991).*
We know that the general Rice formula giving the factorial moments of the number of level crossings of a stochastic process satisfying some conditions can hold whether finite or not.
In the search of conditions for the finiteness of moments of the number of crossings, Nualart and Wschebor (see [113]) proposed some sufficient conditions in the case of a general stochastic process, that reduce in the Gaussian case to:

**Theorem 2.11** *(Nualart and Wschebor, 1991)*
*If $X = (X_t, t \in I \subset I\!\!R)$ is a Gaussian process having $C^\infty$ paths and such that*

$var(X_t) \geq a > 0$, $t \in I$, then $M_t^k(u) < \infty$ for every level $u \in \mathbb{R}$ and every $k \in \mathbb{N}^*$.

- *Azaïs and Wschebor (2001).*
The computation of the factorial moments of the crossings of a process is still a subject of interest. In particular, when expressing these factorial moments by means of Rice integral formulas (of the type of (2.16) in the case of the 2nd factorial moment of the zero crossings, for instance), there arises the problem of describing the behavior of the integrands (appearing in these formulas) near the diagonal; it is still an open (and difficult) problem, even though partial answers have been provided, as for instance by Azaïs and Wschebor (see [10], [12], and references therein). Let us give an example of the type of results they obtained, which helps to improve the efficiency of the numerical computation of the factorial moments, in spite of the restrictive conditions.

**Proposition 2.1** *(Azaïs and Wschebor, 2001)*
*Suppose that $X$ is a centered Gaussian process with $C^{2k-1}$ paths ($k$ integer $\geq 2$), and that for each pairwise distinct values of the parameter $t_1, t_2, \cdots, t_k \in I$ the joint distribution of $(X_{t_h}, X'_{t_h}, \cdots, X_{t_h}^{2k-1}, \ h = 1, 2, \cdots, k)$ is non degenerate. Then, as $t_1, t_2, \cdots, t_k \to t^*$,*

$$\int_{[0,+\infty)^k} x'_1 \cdots x'_k \, p_{X_{t_1}, \cdots, X_{t_k}, X'_{t_1}, \cdots, X'_{t_k}}(0, \cdots, 0, x'_1 \cdots x'_k) dx'_1 \cdots dx'_k \approx J_k(t^*) \prod_{1 \leq i < j \leq k} (t_j - t_i)^4,$$

*where $J_k(.)$ is a continuous non-zero function.*

*2.1.4. Two reference methods in the Gaussian extreme value theory*

- *The normal comparison method.*
The main tool in the Gaussian extreme value theory, and maybe one of the basic important tools of the probability theory, has certainly been the so-called normal comparison technique.
This method, used in the Gaussian case, bounds the difference between two standardized normal distribution functions by a convenient function of their covariances. This idea seems intuitively reasonable since the finite dimensional distributions of a centered stationary Gaussian process is determined by its covariance function.
It was first developed by Plackett in 1954 (see [125]), by Slepian in 1962 (see [146]), then by Berman in 1964 and 1971 (see [22]) and by Cramér in 1967 (see [38]) in the independent or midly dependent cases. An extension of this method to the strongly dependent case was introduced in 1975 by Mittal and Ylvisaker (see [109] and also the 84s paper [108] of Mittal for a review on comparison techniques).

**Theorem 2.12** *(Normal comparison theorem; Slepian, 1962)*
*Let $\xi_1 = (\xi_1(t), t \in T)$ and $\xi_2 = (\xi_2(t), t \in T)$ (with $T$ parameter set) be two separable normal processes (possessing continuous sample functions) having* *same mean and same variance functions, with respective covariance function $\rho_1$ and $\rho_2$. Suppose that $\rho_1(s,t) \leq \rho_2(s,t), \ s,t \in T$.*

*Then for any $x$, $\ \mathbb{P}\left[\sup_{t \in T} \xi_1(t) \leq x\right] \leq \mathbb{P}\left[\sup_{t \in T} \xi_2(t) \leq x\right]$.*

The proof constitutes a basis for proofs of the Berman inequality and a whole line of its generalizations, and is part of what we call the normal comparison method.

Let us illustrate it by considering a pair of Gaussian vectors of dimension $n$, $X_1$ and $X_2$, that we suppose to be independent, with respective distribution functions $F_i$, $i = 1, 2$, density functions $\varphi_i$, $i = 1, 2$, and covariance matrices $\Sigma_i = ((\sigma_i(j,k))_{j,k})$, $i = 1, 2$ such that $\sigma_1(j,j) = \sigma_2(j,j)$. Then the covariance matrix $\Sigma_h := h\Sigma_2 + (1-h)\Sigma_1$ is positive definite. Let $f_h$ and $F_h$ be respectively the n-dimensional normal density and distribution function based on $\Sigma_h$.

Let recall the following equation, discovered by Plackett in 1954 (see [125]), recorded by Slepian in 1962 (see [146]) and proved later in a simpler way than the one of the author, by Berman in 1987 (see [22]):

**Theorem 2.13** *(Plackett partial differential equation, 1954)*
*Let $\Phi_\Sigma$ be the centered normal density function with covariance matrix $\Sigma = (\sigma_{ij})_{i,j}$. Then*

$$\frac{\partial \Phi_\Sigma}{\partial \sigma_{ij}} \ = \ \frac{\partial^2 \Phi_\Sigma}{\partial x_i \partial x_j}, \quad i \neq j.$$

This equation will help to compute the difference between the two normal distribution functions $F_i$, $i = 1, 2$. Indeed, $x$ being in this case a real vector with coordinates $(x_i)_{1 \leq i \leq n}$, we have

$$F_2(x) - F_1(x) = \int_0^1 \frac{d}{dh} F_h(x) dh = \int_0^1 \sum_{i>j} \frac{\partial}{\partial \sigma_h(i,j)} F_h(x) \frac{d\sigma_h(i,j)}{dh} dh \ =$$

$$\sum_{i>j} (\sigma_2(i,j) - \sigma_1(i,j)) \int_0^1 \oint_{-\infty}^{x_k} f_h(y_1, \cdots, y_{i-1}, x_i, y_{i+1}, \cdots, y_{j-1}, x_j, y_{j+1}, \cdots, y_n) \prod_{k \neq i,j} dy_k dh,$$

where $\oint_{-\infty}^{x_k}$ represents the integral of order $n-2$: $\int_{-\infty}^{x_1} \cdots \int_{-\infty}^{x_{i-1}} \int_{-\infty}^{x_{i+1}} \cdots \int_{-\infty}^{x_{j-1}} \int_{-\infty}^{x_{j+1}} \cdots \int_{-\infty}^{x_n}$.

Thus, if for all $i \neq j$, $\sigma_2(i,j) - \sigma_1(i,j) \geq 0$, then $F_2(x) - F_1(x) \geq 0$.

The same results hold in the case of Gaussian separable functions of arbitrary kind.

From those last two results, Berman obtained:

**Theorem 2.14** *(Berman inequality, 1964-1992)*
*Suppose that $(X_{1i}, 1 \leq i \leq n)$ are standard normal random variables with covariance matrix $\Lambda_1 = (\Lambda_{i,j}^1)$ and $(X_{2i}, 1 \leq i \leq n)$ similarly with covariance matrix $\Lambda_2 = (\Lambda_{i,j}^2)$; let $\rho_{ij} = \max(|\Lambda_{i,j}^1|, |\Lambda_{i,j}^2|)$ and $\rho = \max_{i \neq j} \rho_{ij}$. Then, for any* *real numbers $x_1, \cdots, x_n$,*

$$|I\!P[X_{1j} \leq x_j, 1 \leq j \leq n] \; - \; I\!P[X_{2j} \leq x_j, 1 \leq j \leq n]| \quad \leq \quad (2.20)$$

$$\frac{1}{2\pi} \sum_{1 \leq i < j \leq n} \frac{|\Lambda_{i,j}^1 - \Lambda_{i,j}^2|}{\sqrt{(1 - \rho_{ij}^2)}} \exp \left\{ -\frac{x_i^2 + x_j^2}{2(1 + \rho_{ij})} \right\} \quad \leq$$

$$\frac{1}{2\pi \sqrt{1 - \rho^2}} \sum_{1 \leq i < j \leq n} |\Lambda_{i,j}^1 - \Lambda_{i,j}^2| \exp \left\{ -\frac{x_i^2 + x_j^2}{2(1 + \rho)} \right\}.$$

In particular these results hold when choosing one of the two sequences with iid r.v., hence the maximum does behave like that of the associated independent sequence; it helps to prove, under some conditions, results on the maximum and on the point process of exceedances of an adequate level $x_n$ of a stationary normal sequence with correlation function $r$; for instance we obtain that the point process of exceedances converges to a Poisson process under the weak dependence condition $r_n \log n \to 0$ (see [85], chap.4) or to a Cox process under the stronger dependence condition $r_n \log n \to \gamma > 0$, or even to a normal process if $r_n \log n \to \infty$ (see [85], chap.6, or [108])).

There is a discussion in Piterbarg (1988 for the Russian version, 1996 for the English one) (see [122]) about two directions in which the Berman inequality can be generalized, on the one hand on arbitrary events, on the other hand for processes and fields in continuous time. Piterbarg points out that it is not possible to carry the Berman inequality (2.20) over to the processes in continuous time as elegantly as it was done for the Slepian inequality (2.12), but provides a solution for Gaussian stationary processes with smooth enough paths (see theorems $C3$ and $C4$, pp.10-12 in [122]) (and also for smooth enough stationary Gaussian fields).

Finally let us mention the last refinements of the Berman inequality (2.20) given by Li and Shao in 2002 (see [89]) that provide an upper bound in (2.20), cleared of the term $(1 - \rho_{ij}^2)^{-1/2}$:

**Theorem 2.15** *(Li and Shao, 2002)*
*Suppose that $(X_{1i}, 1 \leq i \leq n)$ are standard normal random variables with covariance matrix $\Lambda_1 = (\Lambda_{i,j}^1)$ and $(X_{2i}, 1 \leq i \leq n)$ similarly with covariance matrix $\Lambda_2 = (\Lambda_{i,j}^2)$;*
*let $\rho_{ij} = \max(|\Lambda_{i,j}^1|, |\Lambda_{i,j}^2|)$ and $\rho = \max_{i \neq j} \rho_{ij}$. Then, for any real numbers $x_1, \cdots, x_n$,*

$$I\!P[X_{1j} \leq x_j, 1 \leq j \leq n] \; - \; I\!P[X_{2j} \leq x_j, 1 \leq j \leq n] \quad \leq \quad (2.21)$$

$$\frac{1}{2\pi} \sum_{1 \leq i < j \leq n} \left( \arcsin(\Lambda_{i,j}^1) - \arcsin(\Lambda_{i,j}^2) \right)^+ \exp \left\{ -\frac{x_i^2 + x_j^2}{2(1 + \rho_{ij})} \right\}$$

*Moreover, for $n \geq 3$, for any positive real numbers $x_1, \cdots, x_n$, and when*

*assuming that* $\Lambda_{i,j}^2 \geq \Lambda_{i,j}^1 \geq 0$ *for all* $1 \leq i, j \leq n$, *then*

$$I\!P[X_{1j} \leq x_j, 1 \leq j \leq n] \ \leq \ I\!P[X_{2j} \leq x_j, 1 \leq j \leq n] \ \leq \tag{2.22}$$

$$I\!P[X_{1j} \leq x_j, 1 \leq j \leq n] \exp\left\{\sum_{1 \leq i < j \leq n} \ln\!\left(\frac{\pi - 2\arcsin(\Lambda_{i,j}^1)}{\pi - 2\arcsin(\Lambda_{i,j}^2)}\right) \exp\left\{-\frac{x_i^2 + x_j^2}{2(1 + \Lambda_{i,j}^2)}\right\}\right\}$$

For other precise versions and extensions of this method, we can also refer to e.g. Leadbetter et al. (see [85]), Tong (see [158]), Ledoux and Talagrand (see [87]) and Lifshits (see [90]).

- *The method of moments, also called the Rice method.*
This method, introduced by Rice to estimate the distribution of the maximum of a random signal, consists in using the first two moments of the number of crossings to estimate the probability of exceeding some given level by a trajectory of a (Gaussian) process, as shown in the lemma below (see [122] p.27 and chap.3). In particular it relies on the fact that the event $(X_0 < x, \max_{0 \leq s \leq t} X_s > x)$ implies the event that there is at least one upcrossing: $(N_t^+(x) \geq 1)$, knowing that the probability of more than one up/down-crossing of the level $x$ becomes smaller as the level becomes larger.
This method works only for smooth processes, but can be extended to non-stationary Gaussian processes (see [136] and [137]) and to non-Gaussian processes.
Let $X = (X_s, s \in [0, t])$ be a.s. continuously differentiable with one-dimensional densities bounded. Then

**Lemma 2.4**

$$0 \ \leq \ I\!E[N_t^+(x)] + I\!P[X_0 \geq x] - I\!P[\max_{0 \leq s \leq t} X_s \geq x]$$

$$\leq \ \frac{1}{2}\left(I\!E[N_t^+(x)\left(N_t^+(x) - 1\right)] + I\!E[N_t^-(x)\left(N_t^-(x) - 1\right)]\right) + I\!P[X_0 \geq x, X_t \geq x],$$

*where* $N_t^-(x)$ *denotes the number of downcrossings of level* $x$ *by* $X$ *on* $[0, t]$ *(recall that* $X_s$ *is said to have a downcrossing of* $x$ *at* $s_0 > 0$ *if for some* $\varepsilon > 0$, $X_s \geq x$ *in* $(s_0 - \varepsilon, s_0)$ *and* $X_s \leq x$ *in* $(s_0, s_0 + \varepsilon))$.

Indeed

$$I\!P[\max_{0 \leq s \leq t} X_s \geq x] = I\!P[X_0 \geq x] + I\!P[X_0 < x, \max_{0 \leq s \leq t} X_s > x]$$

and

$$I\!P[X_0 < x, \max_{0 \leq s \leq t} X_s > x] = I\!P[X_0 < x, N_t^+(x) = 1] + I\!P[X_0 < x, N_t^+(x) \geq 2]$$

$$= I\!P[N_t^+(x) = 1] - I\!P[X_0 \geq x, N_t^+(x) = 1] + I\!P[X_0 < x, N_t^+(x) \geq 2]$$

$$= I\!P[N_t^+(x) = 1] + I\!P[N_t^+(x) \geq 2] - I\!P[X_0 \geq x, N_t^+(x) \geq 1],$$

with $\quad I\!\!P[N_t^+(x) = 1] = I\!\!E[N_t^+(x)] - \sum_{k=2}^{\infty} k\, I\!\!P[N_t^+(x) = k]$

and $I\!\!P[X_0 \geq x, N_t^+(x) \geq 1] \leq I\!\!P[X_0 \geq x, X_t \geq x] + I\!\!P[N_t^-(x) \geq 2]$.

Now we can conclude since we also have $I\!\!P[N_t^-(x) \geq 2] \leq I\!\!E[N_t^-(x)\left(N_t^-(x)-1\right)]/2$ and

$$\sum_{k=2}^{\infty} k\, I\!\!P[N_t^+(x) = k] - I\!\!P[N_t^+(x) \geq 2] \leq I\!\!E[N_t^+(x)\left(N_t^+(x)-1\right)]/2. \ \ \square.$$

Estimates previously proposed for the (factorial) moments can then be used at this stage of calculation.

Recently, Azaïs and Wschebor (see [13]) adapted this method to express the distribution of the maximum of a one-parameter stochastic process on a fixed interval (in particular in the Gaussian case) by means of a series (called Rice series) whose terms contain the factorial moments of the number of upcrossings, and which converges for some general classes of Gaussian processes, making the Rice method attractive also for numerical purpose.

### *2.2. Crossings and local time*

#### *2.2.1. Representations of the number of crossings*

Let $X = (X_s, s \geq 0)$ be a stationary Gaussian process defined on a probability space $(\Omega, \mathcal{F}, I\!\!P)$, with mean zero, variance one and correlation function $r$ such that $-r''(0) = 1$. We are first interested in having a representation of the number of $x$-crossings of $X$ as a sum of multiple Wiener-Itô integrals or in terms of Hermite polynomials, where the $n$th Hermite polynomial $H_n$ can be defined as

$$exp(tx - \frac{t^2}{2}) = \sum_{n=0}^{\infty} H_n(x)\frac{t^n}{n!}$$

or as

$$H_n(x) = (-1)^n e^{x^2/2} \frac{d^n}{dx^n}(e^{-x^2/2}), \qquad n \geq 0.$$

Let $W$ be the standard Brownian Motion (or Wiener process).

Let $H(X)$ denote the space of real square integrable functionals of the process $X$.

Recall that $H(X) = \bigoplus_{n=0}^{\infty} \mathcal{H}_n$, $\mathcal{H}_n$ being the Wiener Chaos, i.e. the closed linear subspace of $L^2(\Omega, \mathcal{F}, P)$ generated by the random variables $\{H_n(W(h)), h \in L^2(I\!\!R, dx), ||h|| = 1\}$, where $W(h)$ is the stochastic integral of $h$ with respect to $W$ and $\mathcal{H}_0$ is the set of real constant functions.

We can make use as well of the multiple Wiener-Itô integral $I_n$ defined as in Major (see [96]), since we have

$$H_n(W(h)) = I_n(h^{\otimes n}) \qquad \forall h \in L^2(I\!\!R, dx), \ s.t. \ ||h|| = 1. \qquad (2.23)$$

This integral operator $I_n$ satisfies the multiplication rule, namely:
for $f_p \in L_s^2(I\!\!R^p, m^p)$ and $g_q \in L_s^2(I\!\!R^q, m^q)$, with

$$L_s^2(I\!\!R^n, m^n) := \left\{ f_n \in L^2(I\!\!R^n, \mathcal{B}(I\!\!R^n), m^n) : \right. \tag{2.24}$$
$$\left. f_n(\underline{\lambda}) = \overline{f_n(-\underline{\lambda})}, \ f_n(\underline{\lambda}) = f_n(\lambda_{p(1)}, \cdots, \lambda_{p(n)}), \forall p \in S_n \right\},$$

$m^n$ denoting the product Borel measure on $I\!\!R^n$, $S_n$ the symetric group of permutations of $\{1, \cdots, n\}$ and $\underline{\lambda} = (\lambda_1, \cdots, \lambda_n)$, then

$$I_p(f_p).I_q(g_q) = \sum_{k=0}^{p \wedge q} \frac{k!(p+q-2k)!}{p!q!} C_p^k C_q^k I_{p+q-2k}(f_p \ \hat{\otimes}_k \ g_q),$$

where $f_p \hat{\otimes}_k g_q$ denotes the average overall permutations of $\lambda$-arguments of the function

$$\int_{I\!\!R^k} f_p(\lambda_1, \cdots, \lambda_{p-k}, x_1, \cdots, x_k) g_q(\lambda_{p-k+1}, \cdots, \lambda_{p+q-2k}, -x_1, \cdots, -x_k) m^k(d\underline{x})$$

and $p \wedge q \equiv min(p,q)$.
We can introduce the Sobolev spaces $I\!\!D^{2,\alpha}$ for $\alpha \in I\!\!R$ as in Watanabe (see [161]).
A functional $f \in H(X)$ with the development

$$f = I\!\!E[f] + \sum_{n=1}^{\infty} I_n(f_n) = \sum_{n=0}^{\infty} I_n(f_n)$$

belongs to $I\!\!D^{2,\alpha}$ if and only if

$$||f||_{2,\alpha}^2 = \sum_{n=1}^{\infty} (1+n)^{\alpha} ||f_n||_2^2 = ||(I-L)^{\alpha/2} f||_2^2 < \infty,$$

where $L$ is the operator defined on $H(X)$ by $Lf = \sum_{n=0}^{\infty} -nI_n(f_n)$.

$L$ coincides with the infinitesimal generator of the Ornstein-Uhlenbeck semigroup $\{T_u, u \geq 0\}$ of contraction operators on $H(X)$ defined by $T_u(f) = \sum_{n=1}^{\infty} e^{-nu} I_n(f_n)$, for any $f = \sum_{n=0}^{\infty} I_n(f_n)$.

In particular we have, for $\beta > 0$,

$$(I-L)^{-\beta} = \frac{1}{\Gamma(\beta)} \int_0^{\infty} e^{-u} u^{\beta-1} T_u du. \tag{2.25}$$

(see [161], p.24, or [110], §1.4).

• *Slud (1991, 1994): MWI integral expansion.*
Multiple Wiener-Itô integrals (MWI) may be a tool to represent and to study non-linear functionals of stationary Gaussian processes, as shown by Kalliampur (see [74], chap.8). Slud first applied the stochastic calculus of MWI integral expansions, in 1991 (see [147]) to express the number of crossings of the mean level by a stationary Gaussian process within a fixed time interval $[0, t]$, the motivation being to obtain probabilistic limit theorems for crossing-counts, then in 1994 (see [148] or [149]) to extend his results to $\mathcal{C}^1$-curve crossings.

**Theorem 2.16** *(Slud, 1991 and 1994)*
*Let $X$ be a mean zero, variance one, stationary Gaussian process with continuous spectral measure, and twice-differentiable correlation function $r$.*
*If $var(N_t(\psi)) < \infty$, then*

$$N_t(\psi(y)) = I\!\!E[N_t(\psi(y))] + \sum_{n=1}^{\infty} I_n(F_n) \quad in \ \ L^2(\Omega), \quad where$$

* *in the case of a given level $x$ (i.e. $\psi(y) = x$, $\forall y$), the mean of $N_t(x)$ is given by the Rice formula (2.1) and*

$$F_n(\underline{\lambda}) = \frac{e^{-\frac{x^2}{2}}}{\pi} \int_0^t e^{is(\lambda_1 + \cdots + \lambda_n)} \times \tag{2.26}$$

$$\sum_{l=0}^{[\frac{n}{2}]} (-r''(0))^{\frac{1}{2}-l} H_{n-2l}(x) \frac{(-1)^{l+1} H_{2l}(0)}{(2l)!(2l-1)} \sum_{1 \le m_1 < \cdots < m_{2l} \le n} \lambda_{m_1} \cdots \lambda_{m_{2l}} ds,$$

* *and in the case of a $\mathcal{C}^1$-curve, the mean of $N_t(\psi)$ is given by the generalized Rice formula (2.7) and*

$$F_n(\underline{\lambda}) = \int_0^t e^{is(\lambda_1 + \cdots + \lambda_n)} \frac{e^{-\frac{u^2}{2}}}{\pi} \left( \sqrt{-r''(0)} H_n(u) - \sum_{l=1}^{n} \frac{i^l}{l!} H_{n-l}(u) \right. \tag{2.27}$$

$$\left. \sum_{1 \le m_1 < \cdots < m_l \le n} \lambda_{m_1} \cdots \lambda_{m_l} \int_0^{(-r''(0))^{-\frac{1}{2}}} e^{-\frac{z^2 y^2}{2}} H_l(-zy) y^{l-2} dy \right)_{u=\psi(s), z=\psi'(s)} ds,$$

*where $\underline{\lambda}$ is an $n$-vector of coordinates $\lambda_i$.*

Note that the functional $N_t(\psi)$ may thus be expressed as the integral on $[0, t]$ of $e^{is(\lambda_1 + \cdots + \lambda_n)}$ multiplied by the formal MWI expansion of the form $\sum_n I_n(f_n(., s))$,

w.r.t. Lebesgue measure.

The proof is mainly based on the discrete approximation method of Cramér and Leadbetter (see [38]), already used to obtain the generalized Rice formula (2.7), by introducing the discrete-time number of crossings $N_\psi(1, 2^{-m})$ defined in (2.8), which increases to $N_\psi(1)$. Since by hypothesis in the theorem $var(N_\psi(1)) < \infty$, so is the limiting variance of $N_\psi(1, 2^{-m})$ (via the monotone convergence theorem); then, because of the orthogonal decomposition of $L^2(\Omega)$, the MWI integrands for $N_\psi(1)$ are the $L^2$ (and *a.e.*) limits of the corresponding integrands for $N_\psi(1, 2^{-m})$.

To provide the MWI integrands for $N_\psi(1, 2^{-m})$, some work on indicators of the type $\mathbb{1}_{(X_0 > c)}$ or $\mathbb{1}_{(X_0 > a)} \mathbb{1}_{(X_h > b)}$, is needed, after having noted that
$$\mathbb{1}_{((X_s - a)(X_{s+h} - b) < 0)} = \mathbb{1}_{(X_s \geq a)} + \mathbb{1}_{(X_{s+h} \geq b)} - 2\mathbb{1}_{(X_s \geq a)} \mathbb{1}_{(X_{s+h} \geq b)}.$$

The main new technical tools used for the study of indicators are a generalization of the Hermite polynomial expansion for the bivariate-normal density in (2.28) and the identity (2.29), both enunciated below because of their own interest.

**Lemma 2.5** *(Slud, 1994)*
$\forall x, y \in \mathbb{R}, \ k, m, n \in \mathbb{N}$ and $|t| < 1$,

$$\sum_{j=0}^{\infty} \frac{t^j}{j!} H_{k+j}(x) H_{m+j}(y) = \frac{(-1)^{k+m} e^{\frac{x^2+y^2}{2}}}{\sqrt{1-t^2}} \frac{\partial^{k+m}}{\partial x^k \partial y^m} e^{-\frac{x^2+y^2-2xyt}{2(1-t^2)}} \qquad (2.28)$$

*and*

$$\left(-\frac{\partial}{\partial x} - \frac{\partial}{\partial y}\right)^n e^{-\frac{x^2+y^2-2xyt}{2(1-t^2)}} = \left(\frac{2}{1+t}\right)^{n/2} H_n\left(\frac{x+y}{\sqrt{2(1+t)}}\right) e^{-\frac{x^2+y^2-2xyt}{2(1-t^2)}}. \qquad (2.29)$$

The case of a constant level is simply deduced from the general case.

Note that Slud used a different method in 1991, when considering directly a constant level, based mainly on properties of generalized hypergeometric functions; indeed he expressed $\mathbb{1}_{((X_s - x)(X_{s+h} - x) < 0)}$ as the sum $\mathbb{1}_{(X_s - x < 0)} \mathbb{1}_{(X_{s+h} - x > 0)} + \mathbb{1}_{(X_s - x > 0)} \mathbb{1}_{(X_{s+h} - x < 0)}$, obtained the MWI expansion for the indicator $\mathbb{1}_{(X_s - x > 0)}$ by using first the Hermite polynomial expansion of this indicator, then by studying the asymptotical behavior of hypergeometrical functions. Then he used the Diagram theorem to express the products of expansions as a sum of Wiener-Itô integrals (see [49] or [96], p.42; another version will also be given in terms of Hermite polynomials by Arcones in 1994 (see [4]), and recalled below) and finally he used Fourier transforms for computations.

● *Kratz and León (1997): Hermite polynomial expansion.*
At the end of the 1990s, Kratz and León (see [80]) proposed a new and direct method to obtain, under some assumptions on the spectral moments of the process, the Hermite polynomial expansion of crossings of any level by a stationary Gaussian process:

**Proposition 2.2** *(Kratz and León, 1997)*
*Let $X$ be a mean zero stationary Gaussian process, with variance one, satisfying*

$$-r''(0) = 1 \ \text{and} \ r^{(iv)}(0) < \infty. \qquad (2.30)$$

*Then the following expansion holds in $L^2(\Omega)$*

$$N_t(x) = \sum_{q=0}^{\infty} \sum_{l=0}^{\left[\frac{q}{2}\right]} b_{q-2l}(x) a_{2l} \int_0^t H_{q-2l}(X_s) H_{2l}(\dot{X}_s) ds, \qquad (2.31)$$

*where $b_k(x) := \frac{1}{k!\sqrt{2\pi}}e^{-x^2/2}H_k(x)$.*

This approach is based on an analytical formula involving the "Dirac function", which formally defines the number of crossings $N_t(x)$ as

$$N_t(x) = \int_0^t \delta_x(X_s)|\dot{X}_s|ds \qquad (2.32)$$

and which can be made precise when approximating the Dirac function; it makes then explicit formulas for MWI expansions much easier to obtain than in Slud's papers, mainly because expanding $|\dot{X}_s|$ in Hermite polynomials in $\dot{X}_s$ rather than in $X_s$ quite simplify the calculations (note that at $s$ fixed, $X_s$ and $\dot{X}_s$ are independent).

Note that the condition (2.30), stronger than the Geman condition (2.9), can from now on be replaced by the Geman condition only because of Kratz and León's recent result (theorem 2.7), which makes the chosen method even more attractive .

Moreover this approach is natural in the sense that formally the Dirac function $\delta_x$ has the generalized Hermite expansion $\delta_x(u) = \sum_{k=0}^{\infty} b_k(x)H_k(u)$ with $b_k(x) = \int_{-\infty}^{\infty} \delta_x(y)\frac{1}{k!}H_k(y)\frac{1}{\sqrt{2\pi}}e^{-y^2/2}dy = \frac{1}{k!}H_k(x)\frac{e^{-x^2/2}}{\sqrt{2\pi}}$ ; then (2.32) has the corresponding development given by lemma 2.6 below, made precise by approximating $\delta_x$ by $\Phi'_{\sigma,x} := \varphi_{\sigma,x}$.

**Lemma 2.6** *(Kratz and León, 1997)*
*Let $X$ satisfy the conditions of proposition 2.2.*
*Let $f \in L^2(\phi(x)dx)$ and let $(c_k , k \geq 0)$ be its Hermite coefficients. One has the following expansion*

$$\int_0^t f(X_s)|\dot{X}_s|ds = \sum_{k=0}^{\infty}\sum_{l=0}^{\infty} c_k a_{2l} \int_0^t H_k(X_s)H_{2l}(\dot{X}_s)ds$$

$$= \sum_{q=0}^{\infty}\sum_{l=0}^{[q/2]} c_{q-2l}a_{2l} \int_0^t H_{q-2l}(X_s)H_{2l}(\dot{X}_s)ds,$$

*where $(a_k , k \geq 0)$ are the Hermite coefficients of the function $|.|$, defined by $a_0 = \left(\frac{2}{\pi}\right)^{1/2}$ and $a_{2l} = \left(\frac{2}{\pi}\right)^{1/2}\frac{(-1)^{l+1}}{2^l l!(2l-1)}$ if $l \geq 1$.*

By defining $\zeta^{K,L} = \sum_{k=0}^{K}\sum_{l=0}^{L} c_k a_{2l}\int_0^t H_k(X_s)H_{2l}(\dot{X}_s)ds$, we note that $(\zeta^{K,L})_{K,L}$ is a Cauchy sequence in $L^2(\Omega)$ and we deduce from the Hermite expansions of $|x|$ and $f(x)$ that $\zeta^{K,L}$ converges to $\int_0^t f(X_s)|\dot{X}_s|ds$ in $L^2(\Omega)$. Then to conclude the proof of the lemma, just notice that the second expansion is a consequence of the orthogonality. □

Let us be more explicit about the proof of proposition 2.2, by giving the main steps.
We will apply the previous result (lemma 2.6) to the function

$$f := \varphi_{\sigma,x}(y) = \sum_{k=0}^{\infty} b_k^{\sigma}(x) H_k(y), \quad \text{with } b_k^{\sigma}(x) = \frac{1}{2\pi\sigma k!} \int_{-\infty}^{\infty} e^{-z^2/2} e^{-\frac{(z-x)^2}{2\sigma^2}} H_k(z) dz,$$

to get the Hermite expansion of $N_t^{\sigma}(x) := \int_0^t \varphi_{\sigma,x}(X_s)|\dot{X}_s|ds$, namely, in $L^2(\Omega)$:

$$N_t^{\sigma}(x) = \sum_{q=0}^{\infty} \sum_{l=0}^{[q/2]} b_{q-2l}^{\sigma}(x) a_{2l} \int_0^t H_{q-2l}(X_s) H_{2l}(\dot{X}_s) ds.$$

Let $\hat{N}_t(x) = \sum_{q=0}^{\infty} \sum_{l=0}^{[q/2]} b_{q-2l}(x) a_{2l} \int_0^t H_{q-2l}(X_s) H_{2l}(\dot{X}_s) ds$,

with $b_k(x) := \lim_{\sigma \to 0} b_k^{\sigma}(x) = \frac{1}{k!} \varphi(x) H_k(x)$.
Now we can write

$$I\!E[(N_t(x) - \hat{N}_t(x))^2] \leq 2(I\!E[(N_t(x) - N_t^{\sigma}(x))^2] + I\!E[(N_t^{\sigma}(x) - \hat{N}_t(x))^2])$$

and prove by Fatou's lemma and by Jensen inequality that $E[(N_t(x) - N_t^{\sigma}(x))^2] \to 0$ and via the chaos decomposition that $E[(N_t^{\sigma}(x) - \hat{N}_t(x))^2]$, as $\sigma \to 0$, to conclude to proposition 2.2. □
Now by using the method of regularization of Wschebor, the result of proposition 2.2 can be extended to a larger class of processes:

**Proposition 2.3** (*Kratz and León, 1997*).
*Let $X$ be a mean zero stationary Gaussian process with variance one and satisfying the Geman condition (2.9). In addition we will assume that*

$$|\theta''(t+h) - \theta''(t)| \leq |h| L_1(h), \quad where \quad \theta''(t) := t L(t), \qquad (2.33)$$

*$L_1(h)$ being an even function belonging to $L^1([0,\delta], dx)$.*
*Then the following expansion holds*

$$\frac{N_t(x)}{\sqrt{-r''(0)}} = \sum_{q=0}^{\infty} \sum_{l=0}^{[\frac{q}{2}]} b_{q-2l}(x) a_{2l} \int_0^t H_{q-2l}(X_s) H_{2l}\left(\frac{\dot{X}_s}{\sqrt{-r''(0)}}\right) ds. \qquad (2.34)$$

Indeed, Wschebor's regularization method allows us to drop the strong condition involving the fourth derivative of the covariance function of the process to replace it by a new smoother condition, which could be named "uniform Geman condition", constituted by the conditions (2.9) and (2.33) together, in the following way.

By introducing the regularized process $X^\varepsilon$ defined by

$$X_t^\varepsilon = \frac{1}{\varepsilon} \int_{-\infty}^{\infty} \varphi\left(\frac{t-s}{\varepsilon}\right) X_s ds, \quad \text{where } \varphi \text{ is a } \mathcal{C}^2 \text{ function with support in } [-1,1],$$

we can check that proposition 2.2 applies to the number $N_t^\varepsilon(x)$ of $x$-crossings associated to $X^\varepsilon$, then we can prove that, under the uniform Geman condition, $N_t^\varepsilon(x) \to N_t(x)$ in $L^2$ as $\varepsilon \to 0$, to conclude on proposition 2.3. Note that for this last convergence, an important step is the use of the diagram formula that we already mentioned (Major, 1981) to show that the partial finite developments of

$\dfrac{N_t^\varepsilon(x)}{\sqrt{-r_\varepsilon''(0)}}$ converge to the same developments in the right hand side of (2.34),

from which we deduce that for each fixed $q$

$$\sum_{l=0}^{\left[\frac{q}{2}\right]} b_{q-2l}(x) a_{2l} \int_0^t H_{q-2l}\left(\frac{X_s^\varepsilon}{\sigma_\varepsilon}\right) H_{2l}\left(\frac{\dot{X}_s^\varepsilon}{\sigma_\varepsilon\sqrt{-r''(0)}}\right) ds \to$$

$$\sum_{l=0}^{\left[\frac{q}{2}\right]} b_{q-2l}(x) a_{2l} \int_0^t H_{q-2l}(X_s) H_{2l}\left(\frac{\dot{X}_s}{\sqrt{-r''(0)}}\right) ds, \text{ in probability as } \varepsilon \to 0. \quad \Box$$

To be more explicit, let us give a version of the Major Diagram formula in terms of Hermite polynomials, which provides the expectations of product of Hermite polynomials over a Gaussian vector, a version that can be found in Breuer and Major (see [29]) or in Arcones (see [4]).

We need to introduce some notations.

Let $G := \{(j,l) : 1 \le j \le p, 1 \le l \le l_j\}$ be the diagram of order $(l_1, \cdots, l_p)$, $V(G)$ the set of vertices $(j,l)$ of the diagram $G$, $\Gamma\{l_1, \cdots, l_p\}$ the set of diagrams of order $(l_1, \cdots, l_p)$, $L_j = \{(j,l) : 1 \le l \le l_j\}$ the $j$th level of $G$, $\{((i,l),(j,m)) : 1 \le i < j \le p, 1 \le l \le l_i, 1 \le m \le l_j\}$ the set of edges $((i,l),(j,m))$ s.t. every vertex is on only one edge. The edges connect vertices of different levels. Given an edge $w = ((i,l),(j,m))$, let $d_1(w) = i$, $d_2(w) = j$.

**Theorem 2.17** : *Diagram formula*
*Let $(X_1, \cdots, X_p)$ be a Gaussian vector mean zero, variance 1 and with $\mathbb{E}[X_i X_j] = r(i,j)$, $1 \le i, j \le p$. Then*

$$\mathbb{E}\left[\prod_{i=1}^p H_{l_i}(X_i)\right] = \sum_{G \in \Gamma\{l_1, \cdots, l_p\}} \prod_{w \in V(G)} r(d_1(w), d_2(w)).$$

• Let us go back to the heuristic formula (2.32) for the number of crossings. Suppose that $X = (X_s, s \ge 0)$ is a mean zero stationary Gaussian process with variance one, that the function of covariance $r$ has two derivatives and satisfies (2.30).

Let $0 = \alpha_0 < \alpha_1 < \dots < \alpha_{k-1} < \alpha_k = t$ be the points where the change of sign of the derivative of $X$ occurs. They are in finite number because the process has a finite fourth spectral moment (condition (2.30)).

There is an $x$-crossing of $X$ between $\alpha_i$ and $\alpha_{i+1}$ if

$$|\mathbb{1}_{(X_{\alpha_{i+1}>x})} - \mathbb{1}_{(X_{\alpha_i}>x)}| = 1,$$

which implies that

$$
\begin{aligned}
N_t(x) &= \sum_{i=0}^{k-1} |\mathbb{1}_{(X_{\alpha_{i+1}>x})} - \mathbb{1}_{(X_{\alpha_i}>x)}|, \\
&= \sum_{i=0}^{k-1} |Y_x(X_{\alpha_{i+1}}) - Y_x(X_{\alpha_i})|,
\end{aligned}
$$

where $Y_x(u) = \mathbb{1}_{[x,\infty)}(u)$ is the Heaviside function whose generalized derivative is the "Dirac function" $\delta_x(u) = \infty$ if $u = x$ and 0 if not.

Therefore we can formally write

$$N_t(x) = \sum_{i=0}^{k-1} \int_{\alpha_i}^{\alpha_{i+1}} \delta_x(X_s)|\dot{X}_s|ds,$$

which gives the formula (2.32).

Note also that (2.32) allows to retrieve formally some classical results, as:

- the Rice formula:

$$\mathbb{E}[N_t(x)] = \int_0^t \mathbb{E}[\delta_x(X_s)]\mathbb{E}|\dot{X}_s| \; ds = te^{-x^2/2} \sqrt{-r''(0)} /\pi \; ;$$

- Cabaña's formula:
  if $g$ is a positive function on $\mathbb{R}$, if $G$ is a primitive of $g$, then (see [32])

$$
\begin{aligned}
\int_{\mathbb{R}} N_t(x)g(x)dx = \int_0^t \left[\int_{\mathbb{R}} \delta_{X_s}(x)g(x)dx\right] |\dot{X}_s| \; ds &= \int_0^t g(X_s) |\dot{X}_s| \; ds \\
&= \sum_{i=0}^{k-1} |G(X_{\alpha_{i+1}}) - G(X_{\alpha_i})| \; .
\end{aligned}
$$

- the Kac formula:
  since $\hat{\delta}_0(t) = 1$, by formally applying the Fourier inversion formula, we have $\delta_0(u) = \frac{1}{2\pi} \int_{\mathbb{R}} cos(tu)dt$,

$$N_t(0) = \int_0^t \delta_0(X_s) |\dot{X}_s| \; ds = \frac{1}{2\pi} \int_{\mathbb{R}} \int_0^t cos(\xi X_s) |\dot{X}_s| \; d\xi ds.$$

• *Kratz (2000)*

If we choose for a Gaussian process the Brownian motion $B$, then its number of $x$-crossings is such that either $N_t(x) = \infty$ or $N_t(x) = 0$ a.s., whereas its local time $L_t(x)$ can be defined formally by

$$L_t(x) = \int_0^t \delta(x)(B_s)ds,$$

and rigourously by a Gaussian approximation (see the table below) or by a uniform one, as $L_t(x) = \lim_{\varepsilon \to 0} \int_0^t \frac{1}{2\varepsilon} \mathbb{1}_{[-\varepsilon,\varepsilon]}(B_s)ds$ a.s. (see for instance [55]).

More generally, let us define the local time $L_t(x)$ of a Gaussian process $X$ as the density of the occupation measure of $X$, as for instance in Berman (see [22]). Its construction by limiting processes is based on the sample path properties of $X$, as first introduced by Levy in the case of the Brownian motion (see for instance [55]).

We may then notice that the notions of number of crossings and local time present some "analogies" in their formulas (heuristic and non-heuristic), as it is shown in the following table:

| Local time $L_t^x$ for $\tilde{X} = (\tilde{X}_s, s \geq 0)$ (for which $L_t^x$ does exist) | Crossings number $N_t^x$ for $X = (X_s, s \geq 0)$ (for which $N_t^x$ does exist) |
|---|---|
| Sojourn time of $\tilde{X}$ on a level $x$ in $[0,t]$: $S_x(t) = \int_0^t \mathbb{1}_{(\tilde{X}_s \geq x)}ds := \int_0^t Y_x(\tilde{X}_s)ds$ | |
| Local time = density of $S_x(t)$ Formally $L_t^x = \int_0^t \delta_x(\tilde{X}_s)ds$ | Number of crossings $N_t^x = \int_0^t \delta_x(X_s)\|\dot{X}_s\|ds$ |

Mathematically

the Heaviside function $Y_x$ and the "Dirac function" $\delta_x$ are approximated respectively by the distribution function and the density of a Gaussian r.v. mean $x$ and variance $\sigma^2$, with $\sigma \to 0$.

Moreover, $\forall f \in \mathcal{C}^\infty$,

| $\int_0^t f(\tilde{X}_s)ds = \int_{-\infty}^\infty f(x)L_t^x dx$ | $\int_0^t f(X_s)\|\dot{X}_s\|ds = \int_{-\infty}^\infty f(x)N_t(x)dx$ Banach-Kac (1925-43) formula (see [15] and [73]) |

Because of this type of correspondance with the local time, we became interested in 'validating mathematically' the heuristic formula (2.32) for crossing-counts by looking for the appropiate Sobolev space, as done for Brownian local time by Nualart and Vives (see [111],[112]), and Imkeller et al. (see [69])).

Those authors proved that $L_t^x \in \begin{cases} \mathbb{D}^{2,\alpha}, & 0 < \alpha < 1/2 \\ \mathbb{D}^{p,\alpha}, & \forall p \geq 2, \ 0 < \alpha < 1/2 \end{cases}$, where $\mathbb{D}^{2,\alpha}$ denotes the Sobolev space over canonical Wiener space obtained via completion of the set of polynomials $F$ with respect to the norms $||F||_{2,\alpha} = ||(I-L)^{\alpha/2}F||_2$, $\alpha \in \mathbb{R}$.

In what concerns the number of crossings $N_t(x)$, whereas the distribution $\delta_x(X_s)$ and the random variable $|\dot{X}_s|$ belong to $\mathbb{D}^{2,\alpha}$ for any $\alpha < -1/2$ and for any $\alpha \in \mathbb{R}$ respectively, the integral in the definition of $N_t(x)$ would appear to have a smoothing effect, showing $N_t(x)$ as a random variable which would belong to $\mathbb{D}^{2,\alpha}$ for any $\alpha < 1/4$ (see [78] and [79]).

### 2.2.2. Approximation of the local time by the number of crossings

Several authors took an interest in the problem of approximating the local time of an irregular process $X$ by the number of crossings of a regularization of this process. One classical regularization is the one obtained by convolution, already presented, defined by

$$X_t^\varepsilon = \psi_\varepsilon * X_t, \quad \text{where} \quad \psi_\varepsilon(.) = \frac{1}{\varepsilon}\psi\left(\frac{\cdot}{\varepsilon}\right) \quad \text{and } \psi \text{ some smooth function.}$$
$$(2.35)$$

● *Wschebor (1984, 1992).*

Wschebor considered this problem for a specific Gaussian process, the Brownian motion (in the multiparametric case, but we will only present the one parameter case).

Let $W = (W_t, t \geq 0)$ be a Brownian motion and let $W_\varepsilon$ be the convolution approximation of $W$ defined by $W_\varepsilon := W * \psi_\varepsilon$, where $\psi_\varepsilon$ is defined in (2.35) with $\psi$, a non-negative $\mathcal{C}^\infty$ function with compact support.

Let $N_\varepsilon^x([a,b])$ be the number of crossings of level $x$ by the regularized process $W_\varepsilon$ on the interval $[a,b]$.

Wschebor (see [162]) showed that, for $x \in \mathbb{R}$,

**Theorem 2.18** *(Wschebor, 1984)*

$$\sqrt{\frac{\pi}{2}}\frac{\varepsilon^{1/2}}{||\psi||_2}N_\varepsilon^x([a,b]) \xrightarrow{L^k} L(x,[a,b]), \quad as \quad \varepsilon \to 0, \quad for \ \ k = 1, 2, \cdots \quad (2.36)$$

*where $L(x,[a,b])$ is the Brownian motion local time at level $x$ on $[a,b]$.*

Later (see [164]), he proved the following:

**Theorem 2.19** *(Wschebor, 1992)*
*For any continuous and bounded function $f$, for a.e. fixed $w$,*

$$\frac{1}{||\psi||_2}\sqrt{\frac{\varepsilon\pi}{2}}\int_{-\infty}^{\infty}f(x)N_\varepsilon(x)dx \; - \; \int_{-\infty}^{\infty}f(x)L(x)dx \; \to \; 0, \quad as \quad \varepsilon \to 0,$$

*where $N_\varepsilon(x) := N_\varepsilon^x([0,1])$ and $L(x) := L(x,[0,1])$.*

• *Azaïs, Florens-Zmirou (1987, 1990)*
Azaïs and Florens-Zmirou, in their 87 paper (see [8]), extended Wschebor's result (2.36) to a large class of stationary Gaussian processes, when considering the $L^2$ convergence and the zero crossings. Under some technical conditions on the Gaussian stationary process $X$ and its regularization $X_\varepsilon := X * \psi_\varepsilon$ (requiring a bit more than the non-derivability of $X$ and giving some bounds on the second and fourth spectral moments of $X_\varepsilon$), on its correlation function $r$ (namely $r$ twice differentiable outside a neighborhood of zero, with bounded variation at zero, $r'$ and $r''$ bounded at infinity) and on the convolution kernel $\psi$ (i.e. $\psi(u)$ and $\psi'(u)$ bounded by a constant time $|u|^{-2}$), they proved that

**Theorem 2.20** *(Azaïs and Florens-Zmirou, 1987)*
*If $X$ admits a local time $L(x,[0,T])$ continuous in $x$ at zero, then*

$$\sqrt{\frac{\pi}{2\lambda_{2,\varepsilon}}}N_\varepsilon^0([0,T]) \;\underset{\varepsilon\to 0}{\to}\; L(0,[0,T]) \quad in \quad L^2,$$

*where $\lambda_{2,\varepsilon}$ denotes the second spectral moment of $X_\varepsilon$.*

Note that Azaïs in 1990 (see [6]) considered more general stochastic processes and provided sufficient conditions for $L^2$-convergence of the number of crossings of some smooth approximating process $X_\varepsilon$ of $X$ (which converges in some sense to $X$) to the local time of $X$, after normalization.

• *Berzin, León, Ortega (1992, 1998), Azaïs and Wschebor (1996).*
Whereas Wschebor proved the a.s. convergence of the variable given in theorem 2.19, León and Ortega studied the $L^2$ convergence of this variable slightly modified, with the rates of convergence, first in 1992 (see [88]), then in collaboration with Berzin, in 1998 (see [24]), under slightly different hypotheses. Indeed, for $X = (X(t), t \in [0,1])$ a centered stationary Gaussian process with covariance function $r$ satisfying $r(t) \sim_{t\to 0} 1 - C|t|^{2\alpha}, 0 < \alpha < 1$, they considered the regularized process $X_\varepsilon = \dfrac{\psi_\varepsilon * X}{var(\psi_\varepsilon * X)}$ (where the kernel $\psi_\varepsilon$ approaches the Dirac function as $\varepsilon \to 0$), and proved that *the convergence of*

$$\frac{1}{\varepsilon^{a(\alpha)}}\int_{-\infty}^{\infty}f(x)\left(\frac{N_\varepsilon(x)}{c(\varepsilon)} - L(x)\right)dx, \qquad with \; c(\varepsilon) = \left(\frac{2}{\pi}\frac{var(\dot{X}_\varepsilon(t))}{var(X_\varepsilon(t))}\right)^{1/2},$$

*($f \in L^4(\phi(x)dx)$ satisfying certain regularity conditions), is in $L^2$ or is the weak convergence, depending on the values of $\alpha$ ($0 < \alpha < 1/4$, $1/4 < \alpha < 3/4$ and $3/4 < \alpha < 1$).*

Meantime, Azaïs and Wschebor considered in 1996 (see [10]) the a.s. convergence of $\int_{-\infty}^{\infty} f(x)N_\varepsilon(x)dx$, for any continuous function $f$ and for $X$ belonging to a larger class of Gaussian processes than the one defined by Berzin, León and Ortega.

Note that we did not mention all the papers on this subject (or closely related to it); indeed much work on the approximation of local times and occupation measures has been done also on specific Gaussian processes such as the fractional Brownian motion (we can quote for instance the work by Révész et al. in the 80s on invariance principles for random walks and the approximation of local times and occupation measures, as well as the work by Azaïs and Wschebor (see [11]) on the first order approximation for continuous local martingales).

### *2.3. Asymptotic results*

#### *2.3.1. The law of small numbers and rate of convergence*

● Due to the extension of the Poisson limit theorem from independent to dependent Bernoulli random variables (see [36]), there has been the same extension in extreme value theory, in particular for the study of the number of (up-)crossings. As the level $x = x(t)$ increases with the length $t$ of the time interval, the upcrossings tend to become widely separated in time provided that there is a finite expected number in each interval. Under some suitable mixing property satisfied by the process $X$ to ensure that the occurrences of upcrossings in widely separated intervals can be considered as asymptotically independent events, we get a limiting Poisson distribution for the number of $x$-upcrossings $N_x^+$ by using the Poisson theorem for dependent random variables (see [160], [38], chap.12, [126], [116], [117], [22], and, for a survey and further minor improvement, [85], chap.8 and 9). We will give here a version of such a result, proposed in Leadbetter et al. (see [85], theorem 9.1.2 p.174).

Under the condition that

$$\mu(\mathbf{x}) := I\!E[\mathbf{N}_{\mathbf{x}}^+((\mathbf{0},\mathbf{1}))] < \infty,$$

then $N_x^+(I) < \infty$ a.s. for bounded $I$, and the upcrossings form a stationary point process $N_x^+$ with intensity parameter $\mu = \mu(x)$. The point process of upcrossings has properties analogous to those of exceedances in discrete parameter cases, namely:

**Theorem 2.21** *(Leadbetter et al.)*
*Let us assume that the correlation function $r$ of the centered stationary normal process $X$ satisfies*

$$r(s) = 1 - \frac{\lambda_2}{2}s^2 + o(s^2) \quad as \quad s \to 0 \tag{2.37}$$

*and the Berman condition given by*

$$r(s)\log s \underset{s \to \infty}{\to} 0 \ . \tag{2.38}$$

*Suppose that $x$ and $t$ tend to $\infty$ in a coordinated way such that*

$$I\!\!EN_x^+(t) = t\mu(x) \underset{t\to\infty}{\to} \tau \quad \text{for some constant} \quad \tau \geq 0. \qquad (2.39)$$

*Define a time-normalized point process $N_*^+$ of upcrossings having points at $s/t$ when $X$ has an upcrossing of $x$ at $s$.*
*Then $N_*^+$ converges to a Poisson process with intensity $\tau$ as $t \to \infty$.*

Note that the asymptotic Poisson character of upcrossings also applies to non-differentiable normal processes with covariance functions satisfying

$$r(\tau) = 1 - C|\tau|^\alpha + o(|\tau|^\alpha), \quad \text{as} \quad \tau \to 0 \quad (\text{with } 0 < \alpha < 2, \text{ and } C \text{ positive constant}), \qquad (2.40)$$

if we consider $\varepsilon$-upcrossings, introduced by Pickands (see [116]) and defined below, instead of ordinary upcrossings (for which under (2.40) their mean number of any level per unit time is infinite).
For a given $\varepsilon > 0$, a process $\xi$ is said to have an $\varepsilon$-upcrossing of the level $x$ at $t$ if $\xi(s) \leq x$, $\forall s \in (t - \varepsilon, t)$ and, $\forall \eta > 0$, $\exists s \in (t, t + \eta)$, $\xi(s) > x$;
so an $\varepsilon$-upcrossing is always an upcrossing while the reverse is not true.
The expectation of the number of $\varepsilon$-upcrossings of $x$ by $\xi$ satisfying (2.40) can be evaluated, and it can be proved that asymptotically this mean number is independent of the choice of $\varepsilon$ for a suitably increasing level $x$, which leads, under the Berman condition and (2.40), to the Poisson result of Lindgren et al. (1975) for the time-normalized point process of $\varepsilon$-upcrossings (see [116] and [85], chap.12 for references and more detail).

• Remark: Another notion related to the one of upcrossings of level $x$ by the process $X$ is the time spent over $x$, called also sojourn of $X$ above $x$, defined in the table of section 2.2.1 by $S_x(t) = \int_0^t I\!\!I_{(\tilde{X}_s \geq x)} ds$. An upcrossing of the level $x$ marks the beginning of a sojourn above $x$. If there is a finite number of expected upcrossings in each interval, then the number of sojourns above $x$ is the same as the number of upcrossings. Under some mixing conditions (recalled in the result of Bermangiven below) and the assumption of a finite expected number of upcrossings in each interval, Volkonskii & Rozanov (see [160]), then Cramér & Leadbetter (see [38]) showed, under weaker conditions, by using the reasoning used in the proof of the Poisson limit for the distribution of upcrossings, that the sojourn above $X$ has a compound Poisson limit distribution (as the sum of a Poisson distributed random number of nearly independent r.v. which are the durations of the sojourns). As early as the 70s, Berman (see [22] for the review of these topics) proposed an alternative to discussing upcrossings; he introduced a method based on Hermite polynomial expansions to study the asymptotic form of the sojourns of $X$ above a level $x$, on the one hand for $x \to \infty$ with fixed $t$, on the other hand for $x, t \to \infty$ in a coordinated way. He considered a larger class of processes (with sample functions not necessarily differentiable), allowing a possible infinite expected number of upcrossings in each finite interval, to prove the compound Poisson limit theorem.

More specifically, by using arguments on moments, Berman (see [22], chap.7) proved that:

*if the level function $x = x(t)$ satisfies the asymptotic condition $x(t) \sim \sqrt{2 \log t}$ as $t \to \infty$ and for a covariance function $r$ such that $1 - r$ is a regularly varying function of index $\alpha$ $(0 < \alpha \leq 2)$ when the time tends to 0, then*
*the limiting distribution of $v(t)S_x(t)$ is a compound Poisson distribution,*
*where $v$ is an increasing positive function determined by the asymptotic form of $1 - r(t)$ for $t \to 0$ and the compounding distribution is uniquely defined in terms of the index of regular variation of $1 - r(t)$ for $t \to 0$.*

● For practical use of the asymptotic theory, it is important to know how faithful and accurate these Poisson approximations are, and in view of applications, we must study the rate of convergence carefully. It took some years before getting information about those rates.

Finally, Pickersgill and Piterbarg (see [120]) established rates of order $t^{-\nu}$ for point probabilities in theorem 2.21, but without giving any information on the size of $\nu$.

Kratz and Rootzén (see [84]) proposed, under the condition that $r(.)$ and $r'(.)$ decay at a specified polynomial rate (quite weak assumption, even if more restrictive than $r(s) \log s \underset{s \to \infty}{\to} 0$), bounds for moments of the number of upcrossings and also for the rate of convergence in theorem 2.21 roughly of the order $t^{-\delta}$, with

$$\delta = \frac{1}{2} \wedge \inf_{s \geq 0} \rho(s), \tag{2.41}$$

where $\rho(.)$ has been defined by Piterbarg (see [121]) by

$$\rho(s) = \frac{(1 - r(s))^2}{1 - r(s)^2 + r'(s)|r'(s)|}. \tag{2.42}$$

Their approach proceeds on the one hand by discretization and blocking to return to discrete parameter cases, and on the other hand by combining the normal comparison method and the Stein-Chen method (see [152] and [36]) as developed by Barbour, Holst and Janson (see [61], [16] and also [52]), this last method imposing to measure the rate of convergence with the total variation distance, defined between integer valued r.v.'s. $X$ and $Y$ by

$$d(X, Y) = \frac{1}{2} \sum_k |I\!P(X = k) - I\!P(Y = k)| = \sup_A |I\!P(X \in A) - I\!P(Y \in A)|.$$

More precisely, the authors achieved the following result on what concerns the rate of convergence, throughout assuming for convenience that $x \geq 1$ and $t \geq 1$:

**Theorem 2.22** *(Kratz and Rootzén, 1997)*
*Suppose $\{X(s); s \geq 0\}$ is a continuous stationary normal process which satisfies (2.11), (2.12) with $\gamma = 2$, (2.13) and (2.14) and that $r(s) \geq 0$ for $0 \leq s \leq t$.*

*Then there are constants $K$ and $K'$ which depend on $r(s)$ but not on $x$ or $t$ such that*

$$d(N_x^+(t), \mathcal{P}(t\mu(x))) \leq K \frac{x^{2+2/\alpha}}{t^\delta} \qquad (2.43)$$

*and, for $t \geq t_0 > 1$,*

$$d(N_x^+(t), \mathcal{P}(t\mu(x))) \leq K' \frac{(\log t)^{1+1/\alpha}}{t^\delta}. \qquad (2.44)$$

The constants $K$ in (2.43) and $K'$ in (2.44) are specified in the authors paper (see [84]).

Let us present briefly the method used to prove such results.

The discretisation consists in replacing the continuous process $\{X(s); s \geq 0\}$ by a sampled version $\{X(jq); j = 0, 1, ...\}$, with $q$ roughly chosen equal to $s^{-1/2}$, which makes extremes of the continuous time process sufficiently well approximated by extremes of its sampled version. Nevertheless it raises a new problem: we may count too many exceedances. So a new parameter $\theta$ is introduced to divide the interval $(0, t]$ into blocks $((k-1)\theta, k\theta]$, $k = 1, ..., t/\theta$ with $\theta$ roughly equal to $t^{1/2}$ (this is what blocking method refers to). This choice of $\theta$ makes the blocks long enough to ensure approximate independence of extremes over disjointed blocks.

It leads to considering $W$ the number of blocks with at least one exceedance:

$W := \sum_{k=1}^{t/\theta} I_k$, where $I_k := \mathbb{1}_{(\max_{jq \in ((k-1)\theta, k\theta]} X(jq) > x)}$. Let $\lambda = \mathbb{E}[W]$.

Then by the triangle inequality for the total variation distance $d$, we have

$$d(N_x(t), \mathcal{P}(t\mu)) \leq d(N_x(t), W) + d(W, \mathcal{P}(\lambda)) + d(\mathcal{P}(\lambda), \mathcal{P}(t\mu)), \qquad (2.45)$$

and we estimate the three terms in the righthand side (RHS) of (2.45) separately.

To estimate $d(N_x(t), W)$, we introduce $N_x^{(q)}(\theta)$ the number of exceedances of $x$ by the sampled process $\{X(jq); jq \in (0, t]\}$ upon an interval of length $\theta$, and the essential part concerns the evaluation of the difference between $N_x^{(q)}(\theta)$ and the probability that we have at least one exceedance by the sampled process on the same interval, i.e. $\mathbb{E}[N_x^{(q)}(\theta)] - \mathbb{E}[I_1]$. To evaluate it, we use lemma 2.3.

The last term of the RHS is bounded in the same way as the first one, using $d(\mathcal{P}(\lambda), \mathcal{P}(t\mu)) \leq |\lambda - t\mu|$ (see e.g. [61], p.12), where $\lambda - t\mu = \frac{t}{\theta}(\mathbb{E}[I_1] - \theta\mu)$.

To study the middle term of the RHS, we use the Stein-Chen method described below. It means that we have to evaluate double sums of covariances of indicators that are grouped into blocks, which implies, because of stationarity, the study of $cov(I_1, I_k), k \geq 2$. A version of the normal comparison lemma given in theorem 2.14 (by taking $x = \min(x_i, 1 \leq i \leq n)$) helps to treat the case $k \geq 3$. When $k = 2$, a new parameter $\theta^*$ is introduced to ensure a quasi-independence between blocks. Then classical normal techniques are used.

Concerning references about the approximation speed in related problems, we could mention the works by Borodin and Ibragimov in the 80s and 90s (see [28]), by Jacod (see [72]), by Perera and Wschebor (see [115]), by Berzin and León (see [23]), etc...

- *The Stein-Chen method.*

In 1970, Stein proposed a new method to obtain central limit theorems for dependent r.v. (see [152]). Chen adapted this method in 1975 to the Poisson approximation (see [36]).

Let $(\Omega, \mathcal{B}, P)$ be a probability space, $E$ the expectation under $I\!P$ and $Z : \Omega \to I\!R$ a r.v. such that $E[Z] < \infty$.

The purpose of such a method is to approximate $E[Z]$ in order to approximate the c.d.f. of a real r.v. $W$ defined on $(\Omega, \mathcal{B}, P)$ by $1\!\!1_{(W \le w_0)} := Z$.

The basic approach of the Stein method consists in:

1. determining $Ker E := \{Y : E[Y] = 0\}$;
2. searching in $Ker E$ for $Z - c$ for some constant $c$;
3. concluding $E[Z] \simeq c$.

Given an integer valued r.v. $W$, the problem is then how to determine $Ker E$, i.e. to characterize the set of functions $h : I\!N \to I\!R$ such that $E[h(W)] = 0$. To solve this problem, a method of "exchangeability" is applied. Let us recall the definition of a pair of exchangeable variables.

Let $(\Omega, \mathcal{B}, I\!P)$ and $(\Omega_1, \mathcal{B}_1, I\!P_1)$ be two probability spaces. $(X, X')$ is an "exchangeable" pair of mappings of $(\Omega, \mathcal{B}, I\!P)$ into $(\Omega_1, \mathcal{B}_1, I\!P_1)$ if

$$\begin{cases} I\!P[X \in A] = I\!P_1(A), \ \forall A \in \mathcal{B}_1 \\ I\!P[X \in A \cap X' \in A'] = I\!P[X \in A' \cap X' \in A], \ \forall A, A' \in \mathcal{B}_1. \end{cases}$$

A pair $(X, X')$ of exchangeable variables satisfies the following property:

*If* $\mathcal{F} = \{F : \Omega^2 \to I\!R, F \text{antisymmetric (i.e. } F(x, x') = -F(x', x)); E|F(X, X')| < \infty\}$,

$\mathcal{X} = \{h : \Omega \to I\!R \text{ measurable, } s.t. \ E|h(X)| < \infty\}$,

*and* $T : \mathcal{F} \to \mathcal{X} \ s.t. \ (TF)(X) = E[F(X, X')|X]$,

*then* $(X, X')$ *satisfies* $E \circ T = 0$, *with the composition defined by* $(E \circ T)(F) :=$ $E[TF]$.

To apply this method to the Poisson approximation, let us recall the following characterization for the Poisson law:

*$W$ is an integer valued r.v., Poisson $\mathcal{P}(\lambda)$ distributed, if and only if*

*for any bounded $f$,* $E[\lambda f(W + 1) - W f(W)] = 0$.

Hence if $E[\lambda f(W + 1) - W f(W)] \simeq 0$ for all bounded function $f$ defined on $I\!N$, then $W$ is nearly distributed as $\mathcal{P}(\lambda)$.

This approximative equality is often easier to check than the direct approximation $E[f(W)] \simeq E[f(W)]$ where $Z$ is $\mathcal{P}(\lambda)$ distributed.

Therefore, let us illustrate *the Stein-Chen method when $W$ is the number of occurrences of a large number of independent events*. The advantage of this method is that the dependent case involves only minor transformations.

Let $(X_i)_{1 \leq i \leq n}$ be $n$ independent r.v. with values in $\{0,1\}$, with (for fixed $n$)

$$p_i := I\!\!P[X_i = 1], \ 1 \leq i \leq n, \text{ and } \lambda := \sum_{i=1}^{n} p_i.$$

Let $W := \sum_{i=1}^{n} X_i$. We proceed in four steps:

1. Construction of an exchangeable pair $(W, W')$.
   Let $(X_i^*)_{1 \leq i \leq n}$ be $n$ independent r.v., independent of the $X_i$ and s.t. $X_i =^d X_i^*$, $\forall 1 \leq i \leq n$.
   Let $U$ be a uniformly distributed r.v. in $\{1, \cdots, n\}$, independent of the $X_i$ and $X_i^*$.
   Let $W' := W - X_U + X_U^*$. Then $(W, W')$ is exchangeable.

2. Definition of an antisymmetric function $F$.
   Let $F(X, X') := f(X')\mathbb{1}_{(X'=X+1)} - f(X)\mathbb{1}_{(X=X'+1)}$, with $f : I\!\!N \to I\!\!R$.
   Then we can apply the property of exchangeable mappings to $(W, W')$, namely $I\!\!E\left[I\!\!E^X[F(W, W')]\right] = 0$, i.e.

   $$I\!\!E[\lambda f(W+1) - W f(W)] = I\!\!E\left[(f(W+1) - f(W)) \sum_{j=1}^{n} p_j \mathbb{1}_{(X_j=1)}\right]. \tag{2.46}$$

3. Artificial solution $U_\lambda h$ (particular case of the general method of Stein).
   Let $h : I\!\!N \to I\!\!R$ be bounded and let the function $U_\lambda h$ be defined by

   $$U_\lambda h(w) := -\sum_{j=w}^{\infty} \frac{(w-j)!}{j!} \lambda^{j-w} \left(h(j) - I\!\!E[h(Z)]\right),$$

   with $Z$ Poisson $\mathcal{P}(\lambda)$ distributed r.v., solution of the equation

   $$\lambda f(w+1) - w f(w) = h(w) - I\!\!E[h(Z)].$$

   Taking $f := U_\lambda h$ in (2.46) provides that for all bounded function $h$ defined on $I\!\!N$,

   $$I\!\!E[h(W)] = I\!\!E[h(Z)] + \sum_{j=1}^{n} p_j^2 \, I\!\!E[V_\lambda h(W_j)],$$

   where $W_j := \sum_{j' \neq j} X_{j'}$ and $V_\lambda h(w) := U_\lambda h(w+2) - U_\lambda h(w+1)$.

4. Estimation of $d_{tv}(\mathcal{L}(W), \mathcal{P}(\lambda))$.
   Let us choose $h$ s.t. $h := h_A$, $A \subset I\!\!N$, with $h_A(w) := \mathbb{1}_{(w \in A)}$; a bound can then be deduced for $V_\lambda h$ and we obtain

   $$|I\!\!P[W \in A] - I\!\!P[Z \in A]| \leq \left(1 \wedge \frac{1}{\lambda}\right) \sum_{j=1}^{n} p_j^2$$

which gives the same upper bound for the distance in variation $d_{tv}(\mathcal{L}(W), \mathcal{P}(\lambda))$ defined by $d_{tv}(\mathcal{L}(X), \mathcal{L}(Y)) = d_{tv}(X, Y) := \sup_{A \subset \mathbb{N}} |\mathbb{P}[X \in A] - \mathbb{P}[Y \in A]|.$

Note that Barbour (see [16]) extended the Stein-Chen method by combining it with coupling techniques (see [140] and [141]) to solve the general problem of Poisson approximations for the distribution of a sum of r.v., not necessarily independent, $\{0, 1\}$-valued. In the case of independent indicators, Deheuvels et al. (see [47] and references therein) combined semigroup theory (see [140]) with coupling techniques to obtain results for the Poisson approximation of sums of independent indicators; it has been shown to give sharper results than Barbour's method (see [77]).

### 2.3.2. Central Limit Theorems

● *Malevich (1969), Cuzick (1976).*
In the 70s, some work was done to prove Central Limit Theorems (CLT) for the number of zero crossings $N_t(0)$ as $t \to \infty$, for instance by Malevich (see [98]) and by Cuzick (see [41]).
Cuzick gave conditions on the covariance function of $X$ which ensure on the one hand a mixing condition at infinity for $X$ and on the other hand a local condition for the sample paths of $X$; those conditions are weaker than the ones given by Malevich to prove the same result, although the same type of proofs is used.

**Theorem 2.23** *(Cuzick, 1976)*
*Assume that $r''$ exists and take $r(0) = 1$.*
*If the following conditions are satisfied:*
*i) $r$, $r'' \in L^2$,*
*ii) the Geman condition (2.9) is verified,*
*iii)* $\liminf\limits_{t \to \infty} \dfrac{var(N_t(0))}{t} = \sigma^2 > 0,$
*then*
$$\frac{(N_t(0) - \mathbb{E}[N_t(0)])}{\sqrt{t}} \xrightarrow{d} \mathcal{N}(0, \sigma^2) \quad, \text{ as } t \to \infty,$$

*where* $\sigma^2 = \dfrac{1}{\pi} \left\{ \sqrt{\lambda_2} + \int_0^\infty \left( \dfrac{\mathbb{E}[|X_0' X_s'| \mid X_0 = X_s = 0]}{\sqrt{1 - r^2(s)}} - \mathbb{E}^2 |X'(0)| \right) ds \right\}.$

Remark: A sufficient condition to have the assumption $(iii)$ and which is directly related to the covariance structure of $X$ can be expressed by

$$\int_0^\infty \frac{r'(t)^2}{1 - r^2(t)} dt \ < \ \frac{\pi}{2} \sqrt{\lambda_2}.$$

*The m-dependent method.*
The proof of this theorem relies on what is called the $m$-dependent method, based on an idea of Malevich, which consists in approximating the underlying

Gaussian process $X$ (and its derivatives when they exist) by an $m$-dependent process $X_m$ (i.e. such that $\mathbb{E}[X_m(s)X_m(t)] = 0$ if $|s-t| > m$), in order to apply the CLT for the $m$-dependent process given in Hoeffding and Robbins (see [60]). Note that the $m$-dependent process $X^m$ is obtained directly from the stochastic representation of $X$ by clipping out the portion $(-\infty, -m/2 \,[\cup]\, m/2, \infty)$ of the integration domain and by normalizing the resulting integral.

• *Piterbarg (1978-1996).*

Applying both the method of comparison and the method of discretization, Piterbarg provided a central limit theorem for the number $N_x^+(t)$ of upcrossings at level $x$ by a Gaussian stationary process $X$ (see [118] or [122]), namely

**Theorem 2.24** *(CLT for the number $N_x^+(t)$ of upcrossings, Piterbarg, 1978)*
*Let $X = (X_s, s \in \mathbb{R})$ be a stationary Gaussian process, mean zero, unit variance, with covariance function $r$ satisfying the Geman condition (2.9) and*

$$\int_0^\infty s\left(|r(s)| + |r'(s)| + |r''(s)|\right) ds \ < \infty. \tag{2.47}$$

*Then*

$$var N_x^+(t) = \sigma^2 t(1 + o(1)) \quad as \quad t \to \infty,$$

*where*

$$
\begin{aligned}
\sigma^2 \ &= \ \int_0^\infty \int_0^\infty \int_0^\infty yz\left(\phi_s(x,y,x,z) - \frac{e^{-x^2}\sqrt{-r''(0)}}{2\pi}exp\left\{\frac{y^2+z^2}{2r''(0)}\right\}\right) dydzds \\
&+ \ \frac{\sqrt{-r''(0)}e^{-x^2}}{2\pi} \ > \ 0.
\end{aligned}
$$

*Also, the central limit theorem holds for $N_x^+(t)$.*

For the computation of the variance of $N_x^+(t)$, it suffices to write
$var(N_x^+(t)) = \mathbb{E}[N_x^+(t)(N_x^+(t) - 1)] + \mathbb{E}[N_x^+(t)] - \mathbb{E}[(N_x^+(t))^2]$ , and to use results on moments, which, combined with the condition (2.47), give the convergence of $\sigma^2$.

To obtain the CLT, Piterbarg proceeded both by discretization and by smoothing approximation, since only a discretized in time approximating process does not satisfy the conditions of the CLT in discrete time, even under (2.47) .

He then introduced a smooth enough process $X_\triangle^\delta$ discretized in time, in order to apply to it under some conditions a known result, a CLT for the number of upcrossings for stationary Gaussian process in a discrete time (see for instance [122]). This approximating process $X_\triangle^\delta$ has been defined as $X_\triangle^\delta(s) :=$
$X^\delta(k\triangle) \ + \ (s-k\triangle)\dfrac{X^\delta((k+1)\triangle) - X^\delta(k\triangle)}{\triangle}$,

for $s \in [k\triangle, (k+1)\triangle), k = 0, 1, \cdots, \triangle$ such that $1/\triangle$ is an integer and $\triangle \to 0$,

and where the smoothed process $X^\delta$ is given by $X^\delta(s) \ := \ \dfrac{1}{\delta}\displaystyle\int_0^\delta X(s+v)dv$.

● *Berman (1971-1992).*

Let us return briefly to the notion of sojourn time. When choosing the level function $x = x(t)$ of the asymptotic order of $\sqrt{2 \log t}$ (as $t \to \infty$), we recalled in the section *Law of small numbers*, that the sojourn time above the level $x$ tends to a compound Poisson limit, since the sojourns are relatively infrequent and their contributions are few but individually relatively substantial. In this case, the local behavior of the correlation function $r$ was determining in the form of the limiting distribution. Now when choosing the level function rising at a slower rate (such that $x(t)/\sqrt{2 \log t}$ is bounded away from 1), the sojourns become more frequent and their contributions more uniform in magnitude, implying with the customary (in the application of CLT) normalization, namely $\dfrac{S_t(x) - I\!\!E[S_t(x)]}{\sqrt{var(S_t(x))}}$ (with $I\!\!E[S_t(x)]$ not depending on $r$), a normal limit distribution. Here the local behavior of $r$ does play in the normalization function, but not in the form of the limiting distribution.

*Berman proved this CLT for two different types of mixing conditions, namely when the covariance function $r$ decays sufficiently rapidly to 0 as $t \to \infty$, or conversely at a sufficiently small rate.*

In the case of a rapid rate of decay of $r$, the mixing condition is based not on the function $r$ itself but on the function $b$ that appears in its spectral representation given by

$$r(t) = \int_{-\infty}^{\infty} b(t+s)\bar{b}(s)ds,$$

when supposing that the spectral distribution is absolutely continuous with derivative $f(\lambda)$, and where $b$ is the Fourier transform in the $L^2$-sense of $\sqrt{f(\lambda)}$. This mixing condition is given by $b \in L^1 \cap L^2$, which means in particular that the tail of $b$ is sufficiently small so that $r$ tends to 0 sufficiently rapidly. The proof of the CLT is based on the $m$-dependent method. Indeed, Berman introduced a family of $m$-dependent stationary Gaussian processes $\{X_m(t); -\infty < t < \infty\}$ to approach uniformly in $t$ the original process $X$ in the mean square sense for large $m$. Then he deduced a CLT to the normalized sojourn of $X$ from the CLT of the $m$-dependent process (established by adapting a blocking method used in the proofs of CLT for dependent r.v.).

In the case of slowly decreasing covariances, the proof of the CLT relies on a method specific to Gaussian processes, based on the expansion of $S_t(x)$ in a series of integrated Hermite polynomials and the method of moments. In fact it is a special case of what is known as a non-central limit theorem for Gaussian processes with long-range dependence (see [50], [156] and the next section *Non central limit theorems*).

● *Slud (1991, 1994)*

Introducing the theoretical tool of Multiple Wiener Itô integrals (MWI's) allowed some authors as Taqqu (see [154]), Dobrushin and Major (see [50]), Giraitis and Surgailis (see [56]), Maruyama (see [103]), Chambers and Slud (see [34]), ... to prove general functional central limit theorems (FCLT) (and non central limit theorems) for MWI expansions. In the 1990s, Slud (see [147]) applied a general central limit theorem of Chambers and Slud (see [34]) to provide Cuzick's CLT for the zero crossings $N_t(0)$ by the process $X$, without needing Cuzick's additional assumptions to get a strictly positive limiting variance. Then he generalized the result to constant levels (see [149], theorem 3.1), as follows.

**Theorem 2.25** *(Slud, 1994)*
*Let $x \in \mathbb{R}$ be arbitrary, and let $X = (X_t, t \geq 0)$ be a mean 0, variance 1, stationary Gaussian process with continuous spectral measure $\sigma$ and twice-differentiable correlation function $r$. Suppose $var(N_t(x)) < \infty$. Assume that $\int_{-\infty}^{\infty} r^2(s)ds < \infty$ if $x = 0$ and that $\int_{-\infty}^{\infty} r(s)ds < \infty$ if $x \neq 0$. Then, as $t \to \infty$,*

$$\frac{1}{\sqrt{t}} \left( N_x(t) - e^{-x^2/2} t \sqrt{-r''(0)}/\pi \right) \xrightarrow{d} \mathcal{N}(0, \alpha^2),$$

*where $\alpha^2 > 0$ is given by the expansion*

$$\alpha^2 = \sum_{n=1}^{\infty} \frac{1}{n!} \mathbb{E}_{m^n}[|f_n(\Lambda)|^2 \mid \Lambda_1 + \cdots + \Lambda_n = 0] \int_{-\infty}^{\infty} r^n(s)ds,$$

*with $f_n \in L_s^2(\mathbb{R}^n, m^n)$ (introduced in (2.24)) defined by*

$$f_n(\underline{\lambda}_n) := \frac{e^{-x^2/2}(e^{i \sum_{j=1}^n \lambda_j} - 1)}{\pi \sum_{j=1}^n \lambda_j} \times$$

$$\sum_{k=0}^{[n/2]} \frac{(-r''(0))^{1-2k}}{(2k)!(2k-1)} H_{n-2k}(x) H_{2k}(0) i^{1+2k} \sum_{1 \leq n_1 < \cdots < n_{2k} \leq n} \lambda_{n_1} \cdots \lambda_{n_{2k}}.$$

Note that Slud extended this result to curve crossings in cases where the curve is periodic and the underlying process has rapidly decaying correlations (see [149], theorem 6.4).

• *Kratz and León (2001): a general method.*
In 2001, Kratz and León (see [82]) introduced a general method which can be applied to many different cases, in particular when the dimension of the index set $T$ is bigger than one, to provide a CLT as general as possible for level functionals of Gaussian processes $X = (X_t, t \in T)$. This method is a combination of two approaches, the one developed by the authors in 1997 (see [80]) and one derived from the work of Malevich (see [98]), Cuzick (see [41]) and mainly Berman (see [22]), that consists in approaching the process $X$ by an $m$-dependent process, in order to be able to use well-known results on $m$-dependent processes. Applying this method, a CLT is given for functionals of $(X_t, \dot{X}_t, \ddot{X}_t, t \geq 0)$, which allows in particular to get immediately the CLT for the number of crossings $N_t(x)$ of $X$, given in Slud (see [149], theorem 3.1).

We suppose that the correlation function $r$ of our (stationary Gaussian) process mean zero, variance one, satisfies

$$r \in L^1 \text{ and } r^{(iv)} \in L^2. \tag{2.48}$$

Note that (2.48) implies that $r'$, $r''$ and $r'''$ belong to $L^2$ as well.

Let $Z_s$ be a r.v. independent of $X_s$ and $\dot{X}_s$ (for each $s$ fixed) such that

$$\frac{\ddot{X}_s}{\sqrt{r^{(iv)}(0)}} = \rho_1 X_s + \rho_2 Z_s, \quad \text{with} \quad \rho_1 = \frac{r''(0)}{\sqrt{r^{(iv)}(0)}}, \rho_2 = \sqrt{1 - \rho_1^2}. \tag{2.49}$$

Note that $X_s$ and $\dot{X}_s$ are independent, as well as $\dot{X}_s$ and $\ddot{X}_s$, which ensures the existence of $Z_s$ with the stated properties.

Let $F_t^X$ be the Hermite expansion given in $H(X)$ by

$$F_t^X = \sum_{q=0}^{\infty} \sum_{0 \leq n+m \leq q} d_{qnm} \int_0^t H_n(X_s) H_{q-(n+m)} \left( \frac{\dot{X}_s}{\sqrt{-r''(0)}} \right) H_m(Z_s) ds, \tag{2.50}$$

with $d_{qnm}$ such that $\forall q \geq 0$,

$$\sum_{0 \leq n+m \leq q} d_{qnm}^2 n! m! (q - (n+m))! < C(q), \tag{2.51}$$

with $(C(q))_q$ some bounded sequence.

Besides the property of stationarity of the process and the orthogonality of the chaos, which will be used to simplify the computations whenever possible, let us give the basic points that constitute this general method.

▷ Let $\mathcal{F}_{Q,t}(X)$ be the finite sum deduced from $\mathcal{F}_t(X)$ for $q = 1$ to $Q$, i.e.

$$\mathcal{F}_{Q,t}(X) = \sum_{q=1}^{Q} \sum_{0 \leq n+m \leq q} \frac{d_{qnm}}{\sqrt{t}} \int_0^t I_{qnm}(s) ds, \quad \text{with}$$

$$I_{qnm}(s) = H_n(X_s) H_{q-(n+m)} \left( \frac{\dot{X}_s}{\sqrt{-r''(0)}} \right) H_m(Z_s).$$

▷ First we show that $\mathcal{F}_t(X)$ can be approximated in $L^2$ by the r.v. $\mathcal{F}_{Q,t}(X)$:

$$\lim_{Q \to \infty} \lim_{t \to \infty} I\!\!E[\mathcal{F}_t(X) - \mathcal{F}_{Q,t}(X)]^2 = 0. \tag{2.52}$$

To prove this convergence, we follow the method developed in Kratz & León (see [81], proof of theorem 1) where two results (cited below), one of Taqqu (see [155], lemma 3.2) and the other of Arcones (see [4], lemma 1), helped respectively for the computation of expectations of the form $I\!\!E \left[ \displaystyle\sum_{0 \leq n+m \leq q} d_{qnm} I_{qnm}(0) \right]^2$ and

$$I\!\!E \left[ \sum_{0 \leq n_1+m_1 \leq q} \sum_{0 \leq n_2+m_2 \leq q} d_{qn_1m_1} d_{qn_2m_2} I_{qn_1m_1}(s) I_{qn_2m_2}(s') \right].$$

**Lemma 2.7** *(Arcones inequality, 1994)*
*Let $X = (X_i, 1 \leq i \leq d)$ and $Y = (Y_i, 1 \leq i \leq d)$ be two Gaussian vectors on $\mathbb{R}^d$, mean zero, such that $\mathbb{E}[X_i X_j] = \mathbb{E}[Y_i Y_j] = \delta_{ij}$, $1 \leq i, j \leq d$.*
*Let $r^{(i,j)} := \mathbb{E}[X_i Y_j]$ and $f$ a function on $\mathbb{R}^d$ with finite second moment and Hermite rank $\tau$, $1 \leq \tau < \infty$, w.r.t. $X$. Recall that the Hermite rank of a function $f$ is defined by*

$$rank(f) := \inf \left\{ \tau : \exists l_j \text{ with } \sum_{j=1}^d l_j = \tau \text{ and } \mathbb{E}\left[ (f(X) - \mathbb{E}[f(X)]) \prod_{j=1}^d H_{l_j}(X_j) \right] \neq 0 \right\}.$$

*Suppose $\Psi := \left( \sup_{1 \leq i \leq d} \sum_{j=1}^d |r^{(i,j)}| \right) \vee \left( \sup_{1 \leq j \leq d} \sum_{i=1}^d |r^{(i,j)}| \right) \leq 1$. Then*

$$\mathbb{E}[(f(X) - \mathbb{E}[f(X)])(f(Y) - \mathbb{E}[f(Y)])] \leq \psi^\tau \mathbb{E}\left[ (f(X) - \mathbb{E}[f(X)])^2 \right].$$

In our case, $\psi$ denotes the supremum of the sum of the absolute values of the off-diagonal terms in the column vectors belonging to the covariance matrix for the Gaussian vector

$$\left( X_0, \frac{\dot{X}_0}{\sqrt{-r''(0)}}, Z_0, X_u, \frac{\dot{X}_u}{\sqrt{-r''(0)}}, Z_u \right).$$

Under the conditions on $r$, we have $\int_0^\infty \psi^2(u) du < \infty$.

**Lemma 2.8** *(Taqqu, 1977)*
*Let $p \geq 2$ and $(X_1, \cdots, X_p)$ be standard Gaussian. Then*

$$\mathbb{E}[H_{k_1}(X_1) \cdots H_{k_p}(X_p)] = \begin{cases} \dfrac{k_1! \cdots k_p!}{2^q (q!)} \displaystyle\sum_\Gamma r_{i_1 j_1} \cdots r_{i_q j_q} & \text{if } \sum_{l=1}^p k_l = 2q; \ 0 \leq k_1; \cdots, k_p \leq q \\ \quad 0 \quad otherwise \end{cases}$$

*where $\displaystyle\sum_\Gamma$ is a sum over all indices $i_1, j_1, \cdots, i_q, j_q \in \{1, 2, \cdots, p\}$ such that $i_l \neq j_l, \forall l = 1, \cdots, p$ and there are $k_1$ indices 1, $k_2$ indices 2, $\cdots$, $k_p$ indices p.*

▷ Now by classical tools (the dominated convergence theorem, Fatou lemma, Cauchy-Schwarz inequality, ...) and Arcones inequality, we can get the limit variance $\sigma^2$ of $\mathcal{F}_{Q,t}(X)$ and prove that it is finite:

$$\lim_{Q \to \infty} \lim_{t \to \infty} \mathbb{E}[\mathcal{F}_{Q,t}(X)]^2 := \sigma^2 = \sum_{q=1}^\infty \sigma^2(q) < \infty, \quad \text{where} \quad (2.53)$$

$$\sigma^2(q) = \sum_{0 \leq n_1 + m_1 \leq q} \sum_{0 \leq n_2 + m_2 \leq q} d_{q n_1 m_1} d_{q n_2 m_2} \int_0^\infty \mathbb{E}[I_{q n_1 m_1}(0) I_{q n_2 m_2}(s)] ds.$$

▷ A version of the $m$-dependent method: Berman's method.

This method consists in approaching the process $X$ by an $m$-dependent process, in order to be able to use well-known results on $m$-dependent processes.

Define $X^\varepsilon$ as a $(1/\varepsilon)$-dependent process to approach the process $X$, as follows. Suppose $X$ has a symmetric spectral density $f$; so its Itô-Wiener representation is given by $X_t = \int_{-\infty}^{\infty} e^{i\lambda t}(f(\lambda))^{1/2}dW(\lambda)$, where $W$ is a complex Wiener process defined on $I\!R$.

Let $\psi$ be defined by $\psi(x) = \varphi * \varphi(x)$, where $*$ denotes the convolution and where $\varphi$ is a twice differentiable even function with support contained in $[-1/2, 1/2]$ (so $\psi$ has support in $[-1, 1]$). We can suppose w.l.o.g. that $||\varphi||^2 = 1$.

By using the Fourier inversion formula we can write

$$\psi(x) = \frac{1}{2\pi}\int_{-\infty}^{\infty} e^{-i\lambda x}\hat{\psi}(\lambda)d\lambda = \frac{1}{2\pi}\int_{-\infty}^{\infty} e^{-i\lambda x}|\hat{\varphi}(\lambda)|^2 d\lambda.$$

Let $\hat{\varphi}_\varepsilon(\lambda) = \frac{1}{2\pi\varepsilon}|\hat{\varphi}(\frac{\lambda}{\varepsilon})|^2$, and then introduce

$$X_t^\varepsilon = \int_{-\infty}^{\infty} e^{i\lambda t}(f * \hat{\varphi}_\varepsilon(\lambda))^{\frac{1}{2}}dW(\lambda),$$

(note that $X_t^\varepsilon$ for fixed $t$ is a standard Gaussian r.v.) and its derivatives denoted by $X_t^{(j)\varepsilon}$, $0 \le j \le 2$, such that

$$X_t^{(j)\varepsilon} = \int_{-\infty}^{\infty} e^{i\lambda t}(i\lambda)^j(f * \hat{\varphi}_\varepsilon(\lambda))^{\frac{1}{2}}dW(\lambda).$$

So we can prove that

$$I\!E[\mathcal{F}_{Q,t}(X) - \mathcal{F}_{Q,t}(X^\varepsilon)]^2 \underset{t\to\infty}{\to} 0 \ , \ \text{with } \varepsilon(t) \underset{t\to\infty}{\to} 0, \qquad (2.54)$$

by using Arcones inequality and some results on the correlation between the process (respectively its derivatives) and the $(1/\varepsilon)$-dependent process associated (respectively its derivatives), obtained when working with the spectral representation of the correlation functions, namely

**Proposition 2.4** *(Kratz and León, 2001)*
**(i)** *For all $0 \le j, k \le 2$, $I\!E[X_{t+.}^{(j)\varepsilon} X_t^{(k)\varepsilon}] = (-1)^k r_\varepsilon^{(j+k)}(.)$ converge uniformly over compacts and in $L^2$ as $\varepsilon \to 0$ towards $r_{j,k}(.) := (-1)^k r^{(j+k)}(.)$.*
*For $j = k = 0$, the convergence takes place in $L^1$ as well.*
**(ii)** *For all $0 \le j, k \le 2$, $r_{j,k}^\varepsilon(.) := I\!E[X_{t+.}^{(j)\varepsilon} X_t^{(k)}]$ converge uniformly over compacts and in $L^2$ as $\varepsilon \to 0$ towards $r_{j,k}(.)$.*

▷ Then it is enough to consider the weak convergence of the sequence $\mathcal{F}_{Q,t}(X^\varepsilon)$ towards a Gaussian r.v. as $t \to \infty$ to get the CLT for $F_t^X$.

We can write

$$\mathcal{F}_{Q,t}(X^\varepsilon) = \sum_{q=1}^Q \sum_{0 \leq n+m \leq q} d_{qnm}^\varepsilon \frac{1}{\sqrt{t}} \int_0^t I_{qnm}^\varepsilon(s) ds$$

$$= \frac{1}{\sqrt{t}} \int_0^t \sum_{q=1}^Q \sum_{0 \leq n+m \leq q} d_{qnm}^\varepsilon I_{qnm}^\varepsilon(s) ds := \frac{1}{\sqrt{t}} \int_0^t f_Q \left( X_s^\varepsilon, \frac{\dot{X}_s^\varepsilon}{\sqrt{-r_\varepsilon''(0)}}, Z_s^\varepsilon \right) ds$$

$$= \frac{1}{\sqrt{t}} \int_0^t \theta_s \left[ f_Q \left( X_0^\varepsilon, \frac{\dot{X}_0^\varepsilon}{\sqrt{-r_\varepsilon''(0)}}, Z_0^\varepsilon \right) \right] ds,$$

where $\theta$ is the shift operator associated to the process and $f_Q(x_1, x_2, x_3) = \sum_{q=1}^Q \sum_{0 \leq n+m \leq q} d_{qnm}^\varepsilon H_n(x_1) H_{q-(n+m)}(x_2) H_m(x_3).$

Hence the weak convergence of $\mathcal{F}_Q(X^\varepsilon)$ towards a Gaussian r.v. is a direct consequence of the CLT for sums of $m$-dependent r.v. (see [60] and [21]), which, combined with (2.52), (2.53) and (2.54), provides the CLT for $\mathcal{F}_t(X)$, namely

**Theorem 2.26** *(Kratz and León, 2001)*
*Under the above conditions, we have*

$$\mathcal{F}_t(X) := \frac{F_t^X - \mathbb{E}[F_t^X]}{\sqrt{t}} \longrightarrow \mathcal{N}(0, \sigma^2) \ as \ t \to \infty,$$

*where* $\mathbb{E}[F_t^X] = t d_{000}$ *and* $\sigma^2 = \sum_{q=1}^\infty \sigma^2(q) < \infty,$ *with*

$$\sigma^2(q) = \sum_{0 \leq n_1 + m_1 \leq q} \sum_{0 \leq n_2 + m_2 \leq q} d_{qn_1 m_1} d_{qn_2 m_2} \int_0^\infty \mathbb{E}[I_{qn_1 m_1}(0) I_{qn_2 m_2}(s)] ds,$$

*and* $\quad I_{qnm}(s) = H_n(X_s) H_{q-(n+m)} \left( \frac{\dot{X}_s}{\sqrt{-r''(0)}} \right) H_m(Z_s).$

**Remark.** Condition (2.48) can of course be weakened to:

$$r \in L^1 \quad \text{and} \quad r'' \in L^2 \tag{2.55}$$

when considering the process $X$ and its first derivative only, as for the number of crossings,
and to

$$r \in L^1 \tag{2.56}$$

when considering the process $X$ only.

• As an application of theorem 2.26, we get back, under the condition (2.55), Slud's CLT for the number of crossings of $X$ already enunciated in theorem 2.25.

Indeed, let us consider the Hermite expansion of $N_t(x)$ given in (2.34), which corresponds to $F_t^X$ with $m := 0$, $d_{qn} := d_{qn0} = b_{q-2l}(x)a_{2l}1_{n=2l}$, $\forall q, n, l \in I\!\!N$ and with $I\!\!E[F_t^X] = td_{000} = tb_0(x)a_0 = \dfrac{t}{\pi}e^{-x^2/2}$. We have $\forall x \in I\!\!R$,

$$\sum_{l=0}^{[q/2]} b_{q-2l}^2(x)a_{2l}^2(2l)!(q-2l)! < C,$$ with $C$ some constant independent of $q$, which is a consequence of proposition 3 in Imkeller et al. (see [69]) given by

**Proposition 2.5** *(Imkeller, Perez-Abreu and Vives, 1995)*
*Let $1/4 \le \alpha \le 1/2$. Then there exists a constant $c$ such that for any $n \in I\!\!N$,*

$$\sup_{x \in I\!\!R} |H_n(x)e^{-\alpha x^2}| \le cn^{-(8\alpha-1)/12}.$$

Therefore we can apply theorem 2.26 to obtain

$$\sqrt{t}\left(\frac{N_t(x)}{t\sqrt{-r''(0)}} - \frac{e^{-x^2/2}}{\pi}\right) \underset{t\to\infty}{\longrightarrow} \mathcal{N}(0, \sigma^2),$$

which can be written as

$$\sqrt{\mathbf{t}}\left(\frac{\mathbf{N_t(x)}}{\mathbf{t}} - \frac{\sqrt{-\mathbf{r''(0)}}}{\pi}\mathbf{e^{-x^2/2}}\right) \underset{\mathbf{t\to\infty}}{\longrightarrow} \mathcal{N}(\mathbf{0}, -\mathbf{r''(0)}\sigma^2),$$

where $\sigma^2$ $(0 < \sigma^2 < \infty)$ is given by $\sigma^2 = \sum_{q=1}^{\infty} \sigma^2(q)$, with

$$\sigma^2(q) = \sum_{n_1=0}^{[q/2]}\sum_{n_2=0}^{[q/2]} b_{q-2n_1}(x)a_{2n_1}b_{q-2n_2}(x)a_{2n_2} \times$$

$$\int_0^\infty I\!\!E\left[H_{2n_1}(X_0)H_{q-2n_1}\left(\frac{\dot{X}_0}{\sqrt{-r''(0)}}\right)H_{2n_2}(X_s)H_{q-2n_2}\left(\frac{\dot{X}_s}{\sqrt{-r''(0)}}\right)\right]ds.$$

We can easily check that $\sigma^2 \ne 0$, by proving that $\sigma^2(2) \ne 0$ (showing that the determinant of the matrix (expressed in the domain of frequencies) associated to $\sigma^2(2)$ is strictly positive). We conclude to Slud's CLT.

### 2.3.3. Non central limit theorems

From the CLT given in the previous section, or even more generally from the literature of the CLT for non-linear functionals of Gaussian processes, we can ask what happens when some of the conditions of those theorems are violated, in particular when we are under a condition of regular long-range dependence. This problem interested many authors, among whom we can cite in chronological order Taqqu (see [154], [156]), Rosenblatt (see [133], [134], [135]), Dobrushin and Major (see [50]), Major (see [97]), Giraitis and Surgailis (see [56]), Ho and Sun (see [63]), Slud (see [148] or [149]). MWI's expansions proved to be quite useful in defining the limiting behavior in non central limit theorems for functionals of Gaussian processes with regular long-range dependence. It is what allows Slud to prove, by using techniques proper to MWI's, in particular Major's non central limit theorem for stationary Gaussian fields with regular long-range dependence (see [96], theorem 8.2), the following non CLT for level-crossing counts.

**Theorem 2.27** *(Slud, 1994)*
*Let X be a stationary Gaussian process, mean 0, variance 1, with continuous spectrum and twice differentiable correlation function with regular long range dependence, i.e.* $r(s) = (1 + |s|)^{-\alpha} L(s)$, *where L is slowly varying at $\infty$ and $0 < \alpha < 1/2$. Assume also that for some $\delta \in (-\infty, \alpha)$ and constant $C < \infty$, for $k \geq 1, 2$ and all $x \geq 0$,*

$$\frac{1}{|r(x)|} \left| \frac{d^k}{dx^k} r(x) \right| \leq C(1 + |x|)^{\delta}.$$

*i) Let $c \neq 0$ be fixed arbitrarily. Then*

$$\frac{N_t(c) - e^{-c^2/2} t \sqrt{-r''(0)}/\pi}{\sqrt{t^{2-\alpha} L(t)}} \quad \xrightarrow{d} \quad \mathcal{N}\left(0, \frac{2e^{-c^2} c^2 (-r''(0))}{(1-\alpha)(2-\alpha)\pi^2}\right), \text{ as } t \to \infty.$$

*ii) When $c = 0$,*

$$\frac{t^{\alpha-1}}{L(t)} \left(N_t(0) - t\sqrt{-r''(0)}/\pi\right) \quad \xrightarrow{d} \quad \frac{-r''(0)}{\pi} \tilde{I}_2 \left(\frac{e^{i(\lambda_1+\lambda_2)} - 1}{i(\lambda_1 + \lambda_2)}\right), \text{ as } t \to \infty,$$

*where $\tilde{I}_2$ denotes the second order MWI integral operator for a stationary Gaussian process $\tilde{X}$ with correlation function $r_0(s) := \int e^{isx} \sigma_0(dx)$ uniquely determined by*

$$\int e^{isx} \frac{(1-\cos x)^2}{x^2} \sigma_0(dx) = \int_0^1 (1-x)|x+s|^{-\alpha} dx, \quad s > 0.$$

## 3. Extensions

The aim of this section is to present examples of possible applications of some of the methods reviewed in the previous section.

### 3.1. CLT for other non-linear functionals of stationary Gaussian processes

As an application of the heuristic of the previous section, in particular of theorem 2.26, we may look at various functionals related to crossing functionals of a stationary Gaussian process $X$, as for instance the sojourn time of $X$ in some interval, the local time of $X$ when it exists, or the number of maxima of the process in an interval.

• *Time occupation functionals.*

The simplest application of theorem 2.26 is when the integrand appearing in $F_t^X$, defined in (2.50), depends only on one variable, without needing other conditions than the smooth one $r \in L^1$ (see the remark following the theorem).

▷ A first example of this type is when $F_t^X$ represents the sojourn of the process $X$ above a level $x$ in an interval $[0, t]$, i.e. when $F_t^X := S_x(t)$.

It is easy to obtain the Hermite expansion of $S_x(t)$ as

$$S_x(t) = \sum_{q=0}^{\infty} d_q \int_0^t H_q(X_s) ds,$$

with $d_0 = 1 - \Phi(x)$ and $d_q = \dfrac{1}{q!} \displaystyle\int_x^{+\infty} H_q(u)\phi(u)du = -\dfrac{1}{q!}H_{q-1}(x)\phi(x), \; \forall q \geq 1$.

Then an application of theorem 2.26 yields the CLT for the sojourn time under the condition $r \in L^1$, namely

$$\frac{S_x(t) - td_0}{\sqrt{t}} \underset{t \to \infty}{\longrightarrow} \mathcal{N}(0, \sigma^2) \; ,$$

with $\sigma^2 = \displaystyle\sum_{q=1}^{\infty} d_q^2 \int_0^{\infty} I\!E[H_q(X_0)H_q(X_s)]ds = \sum_{q=1}^{\infty} q! d_q^2 \int_0^{\infty} r^q(s)ds$.

▷ Another example already discussed is the local time $L_t^x$ for $X$ in the level $x$ (when it exists) (see [22]).

Its Hermite expansion is given by

$$L_t^x = \phi(x) \sum_{k=0}^{\infty} \frac{H_k(x)}{k!} \int_0^t H_k(X_s) ds = \sum_{k=0}^{\infty} l_k \int_0^t H_k(X_s) ds.$$

and again theorem 2.26 allows to retrieve its asymptotical normal behavior under the condition $r \in L^1$, namely

$$\frac{L_t^x - t\phi(x)}{\sqrt{t}} \underset{t \to \infty}{\longrightarrow} \mathcal{N}(0, \tilde{\sigma}^2) \; ,$$

with $\tilde{\sigma}^2 = \displaystyle\sum_{k=1}^{\infty} l_k^2 \int_0^{\infty} I\!E[H_k(X_0)H_k(X_s)]ds = \phi^2(x) \sum_{k=1}^{\infty} \frac{1}{k!} \int_0^{\infty} r^k(s)ds$.

• *Number of maxima in an interval.*

One of the main concerns of extreme value theory is the study of the maximum $\left(\underset{t \in [0,T]}{\max} X_t\right)$ of a real-valued stochastic process $X = (X_t, t \in [0, T])$ having continuous paths, in particular the study of its distribution $F$.

There is an extensive literature on this subject, going mainly in three directions, according to Azaïs and Wschebor (see [12]): one looking for general inequalities for the distribution $F$, the other describing the behavior of $F$ under various asymptotics, and the last one studying the regularity of the distribution $F$. For more references, see [12].

In the discussion of maxima of continuous processes, the upcrossings of a level play an important role, as was the case in the discrete case between maxima of sequences and exceedances of level $u_n$ through the equivalence of the events $(M_n^{(k)} \leq u_n) = (N_n^+ < k)$, $M_n^{(k)}$ being the $k$th largest value of the r.v. $X_1, \cdots, X_n$ and $N_n^+$ the point process of exceedances on $(0, 1]$ (i.e. $N_n^+ = \#\{i/n \in (0, 1] : X_i > u_n, 1 \leq i \leq n\}$).

In the continuous case, we have already seen that crossings and maxima are closely related when describing the Rice method, in particular with lemma 2.4 providing bounds of $I\!P[\max_{0 \leq s \leq t} X_s \geq x]$ in terms of factorial moments of the level $x$-crossings of the process $X$.

Note that Cramér (1965) noted the connection between $u$- upcrossings by $X$ and its maximum, e.g. by $\{N_u(T) = 0\} = \{M(T) \leq u\} \cup \{N_u(T) = 0, X(0) > u\}$, which led to the determination of the asymptotic distribution of the maximum $M(T)$.

Recently, when looking for a reasonable way based upon natural parameters of the considered process $X$ to compute the distribution of the maximum of $X$, Azaïs and Wschebor (see [13]) established a method, based upon expressing the distribution of the maximum $I\!P[\max_{0 \leq s \leq 1} X_s > x]$ of the process satisfying some regularity conditions, by means of Rice series, whose main $k$th term is given by $(-1)^{k+1}\nu_k/k!$, $\nu_k$ denoting the $k$th factorial moment of the number of upcrossings. This method, named 'the Rice method revisited' because inspired by previous works such as the one of Miroshin (see [106]), can be applied to a large class of processes, and allows a numerical computation of the distribution in Gaussian cases more powerful in many respects than the widely used Monte-Carlo method, based on the simulation of the paths of the continuous parameter process.

Note also another useful connection with the maximum $\left(\max_{0 \leq s \leq t} X_s\right)$ of a process $X$, which is the sojourn time. Indeed, we can write

$$\left(\max_{0 \leq s \leq t} X_s \leq x\right) \quad \Leftrightarrow \quad (S_x(t) = 0).$$

Berman uses this equivalence between the events $\left(\max_{0 \leq s \leq t} X_s > x\right)$ and $(S_x(t) > 0)$ to study the maximum of the process $X$ (see [22], chap.10).

Here we are interested in the number of local maxima of a stationary Gaussian process lying in some interval, and more specifically in its asymptotical behavior, one of the motivations being the applications to hydroscience.

More precisely let $M^X_{[\beta_1,\beta_2]}$ be the number of local maxima of $X_s, 0 \le s \le t$, lying in the real interval $[\beta_1, \beta_2]$ and let $r(0) = 1$.

Kratz & León (1997) provided under the condition $r^{(vi)}(0) < \infty$ the Hermite expansion of $M^X_{[\beta_1,\beta_2]}$, by adapting the proof of proposition 2.2 on the number of crossings, and with the change of variables (2.49).

Formally, $M^X_{[\beta_1,\beta_2]} = \int_0^t 1_{[\beta_1,\beta_2]}(X_s)\delta_0(\dot{X}_s) \mid \ddot{X}_s \mid 1_{[0,\infty)}(\ddot{X}_s)ds$; more precisely:

**Theorem 3.1** *(Kratz and León, 1997)*
*Under the condition $-r^{(vi)}(0) < \infty$, we have*

$$
M^X_{[\beta_1,\beta_2]} = -\sqrt{\frac{r^{(iv)}(0)}{-r''(0)}} \sum_{q=0}^\infty \sum_{0 \le n+m \le q} \delta_{nm} \frac{H_{q-(m+n)}(0)}{(q-(m+n))!\sqrt{2\pi}}
$$

$$
\int_0^t H_n(X_s)H_{q-(n+m)}\left(\frac{\dot{X}_s}{\sqrt{-r''(0)}}\right)H_m(Z_s)ds
$$

*where $\delta_{nm}$ is defined by*

$$
\delta_{nm} = \frac{1}{n!m!}\int_{\beta_1}^{\beta_2}\!\!\int_{I\!R}(\rho_1 x + \rho_2 z)1_{(-\infty,0)}(\rho_1 x + \rho_2 z)H_n(x)H_m(z)\phi(x)\phi(z)dxdz, \quad (3.57)
$$

*$\rho_1$ and $\rho_2$ satisfying (2.49).*

As for the number of crossings, the condition appearing in the theorem could be weakened by taking a similar condition to (2.33) for the fourth derivative of $r(.)$; it would then provide that $I\!E(M^X_{[\beta_1,\beta_2]})^2 < \infty$.

Previously two cases of application of theorem 2.26 were considered, on the one hand when the integrand appearing in $F^X_t$, defined in (2.50), depends only on one variable (under the simple condition $r \in L^1$), on the other hand when the integrand depends on two variables (under the conditions $r \in L^1$ and $r'' \in L^2$), case of our main study, the number of crossings, but which also concerns any convex combination of it as, for instance, Cabaña estimator (see [32]) of the second spectral moment (slightly modified) given by $\gamma = \frac{\pi}{t}\int_{-\infty}^\infty N_t(x)d\alpha(x)$ (with $\alpha(.)$ a distribution function on $I\!R$), and which has been studied by Kratz & León (see [81]).

To apply theorem 2.26 to obtain a CLT for the number of local maxima $M^X_{[\beta_1,\beta_2]}$ made us consider the last possible case, i.e. when the integrand appearing in $F^X_t$ depends on three variables.

Note that the condition $r^{(vi)}(0) < \infty$ and $r \in L^1$ imply the condition (2.48) of theorem 2.26 and that we can write $M^X_{[\beta_1,\beta_2]} = F^X_t$ when taking

$d_{qnm} = -\sqrt{\dfrac{r^{(iv)}(0)}{-r''(0)}}\dfrac{\delta_{nm}H_{q-(m+n)}(0)}{(q-(m+n))\sqrt{2\pi}}$. Moreover we can easily check that
$\sum_{q=0}^\infty d^2_{qnm}n!m!(q-(n+m))! < C$, with $C$ some constant, by using again proposition 2.5.

Therefore theorem 2.26 yields a CLT for $M^X_{[\beta_1,\beta_2]}$, result obtained by Kratz & León in 2000 (see [81]) with another method which consisted mainly in adapting and verifying conditions (a1)-(a3) in the 1994 Slud's paper (see [149], p.1362). But, as already stated by Berman (see [21], pp.62-63), it wouldn't have been possible to benefit of this method when working in a dimension higher than one, which confirmed us in the idea of finding a more general method as the one put forward to demonstrate theorem 2.26.

**Theorem 3.2** *(Kratz and León, 2000)*
*Under the conditions $-r^{(vi)}(0) < \infty$ and $r \in L_1$,*

$$\frac{M^X_{[\beta_1,\beta_2]} - I\!\!E[M^X_{[\beta_1,\beta_2]}]}{\sqrt{t}} \xrightarrow{d} N(0,\sigma^2) \text{ as } t \longrightarrow \infty,$$

*where*

$$
\begin{aligned}
I\!\!E\left(M^X_{[\beta_1,\beta_2]}\right) &= \frac{t}{\sqrt{2\pi}} \frac{\sqrt{r^{(iv)}(0)}}{\sqrt{-r''(0)}} \left[ \rho_1 \left( \phi(\beta_2)\Phi\left(-\frac{\rho_1}{\rho_2}\beta_2\right) - \phi(\beta_1)\Phi\left(-\frac{\rho_1}{\rho_2}\beta_1\right) \right) \right. \\
&\quad + \left. \frac{1}{\sqrt{2\pi}}\sqrt{\rho_2^2 + \rho_1^2} \left\{ \Phi\left(\beta_2\sqrt{1+\frac{\rho_1^2}{\rho_2^2}}\right) - \Phi\left(\beta_1\sqrt{1+\frac{\rho_1^2}{\rho_2^2}}\right) \right\} \right], \quad (3.58)
\end{aligned}
$$

$\rho_i$ $(i=1,2)$ *being defined in (2.49), and $\sigma^2$ is given by*

$$
\sigma^2 = \frac{r^{(iv)}(0)}{-r''(0)2\pi} \sum_{q=1}^{\infty} \sum_{0 \le n_1 + m_1 \le q} \sum_{0 \le n_2 + m_2 \le q} \delta_{n_1 m_1} \delta_{n_2 m_2} \frac{H_{q-(m_1+n_1)}(0) H_{q-(m_2+n_2)}(0)}{(q-(m_1+n_1))!(q-(m_2+n_2))!} \times
$$

$$
\int_0^{+\infty} I\!\!E\left[ H_{n_1}(X_0) H_{q-(n_1+m_1)}\left(\frac{\dot{X}_0}{\sqrt{-r''(0)}}\right) H_{m_1}(Z_0) H_{n_2}(X_s) H_{q-(n_2+m_2)}\left(\frac{\dot{X}_s}{\sqrt{-r''(0)}}\right) H_{m_2}(Z_s) \right] ds.
$$

This last result has applications to hydroscience, in particular in the study of random seas.

Ocean waves have been studied for many years by many authors, in particular through a stochastic approach by M. Ochi (see [114], for references also). In this approach, the randomly changing waves are considered as a stochastic process (so that it is possible to evaluate the statistical properties of waves through the frequency and probability domains) and in deep water as a Gaussian random process. This Gaussian property was first found by Rudnick in 1951 through analysis of measured data obtained in the Pacific Ocean. This process is assumed to be stationary. It is a common hypothesis for the study of ocean waves since the sea does behave in a stationary way when observed over short periods of time (see [93], [114]).

From the general CLT for the number of maxima of a stationary Gaussian process given in theorem 3.2, we can deduce the asymptotical behavior of functionals related to ocean waves, such as the number of local positive maxima of waves (studied when having a non-narrow-band-spectrum random process (see [114], §3.3)), or wave amplitude, or even the amplitude associated with acceleration (see [114], §4.3)). As an example, let us give the asymptotic behavior of the number of local positive maxima of waves.

Indeed, as shown by Ochi, consideration of the negative maxima is unnecessary in the maximum-point-process of waves for engineering problems, because the negative maxima do not contribute to the largest peak value (the extreme wave amplitude) expected to occur in a certain number of observations.

Let us give the asymptotic behavior of the ratio $\dfrac{M^X_{[\beta,+\infty)}}{M^X_{[0,+\infty)}}$, with $\beta \geq 0$.

**Corollary 3.1** *(Number of local positive maxima of waves, Kratz and León, 2000)*

$$\sqrt{T}\left(\frac{M^X_{[\beta,+\infty)}}{M^X_{[0,+\infty)}} - \frac{I\!E[M^X_{[\beta,+\infty)}]}{I\!E[M^X_{[0,+\infty)}]}\right) \xrightarrow{d} N(0,\sigma^2) \ as \ T \longrightarrow \infty,$$

*where $\sigma^2$ is given by*

$$\sigma^2 = \frac{r^{(iv)}(0)}{-r''(0)2\pi C_0^4} \sum_{q=1}^{\infty} \sum_{0 \leq n_1 + m_1 \leq q} \sum_{0 \leq n_2 + m_2 \leq q} (C_0 \delta_{n_1 m_1}(\beta) - C_\beta \delta_{n_1 m_1}(0)) \ \times$$

$$(C_0 \delta_{n_2 m_2}(\beta) - C_\beta \delta_{n_2 m_2}(0)) \frac{H_{q-(m_1+n_1)}(0) H_{q-(m_2+n_2)}(0)}{(q-(m_1+n_1))!(q-(m_2+n_2))!} \ \times$$

$$\int_0^{+\infty} I\!E\left[H_{n_1}(X_0)H_{q-(n_1+m_1)}\left(\frac{\dot{X}_0}{\sqrt{-r''(0)}}\right)H_{m_1}(Z_0)H_{n_2}(X_s)H_{q-(n_2+m_2)}\left(\frac{\dot{X}_s}{\sqrt{-r''(0)}}\right)H_{m_2}(Z_s)\right] ds,$$

*with $C_y := I\!E[M^X(y)] = \dfrac{\sqrt{r^{(iv)}(0)}}{2\pi\sqrt{-r''(0)}}\left(1 - \Phi\left(y\sqrt{1 + \dfrac{\rho_1^2}{\rho_2^2}}\right)\right) + \sqrt{\dfrac{-r''(0)}{2\pi}}\phi(y)\Phi\left(-\dfrac{\rho_1}{\rho_2}y\right),$*

*and in particular $C_0 = \dfrac{\sqrt{r^{(iv)}(0)} - r''(0)}{4\pi\sqrt{-r''(0)}}.$*

Remark: for statistical purpose, a consistent estimator of the asymptotic variance has also been proposed (see [81]).

### 3.2. Multidimensional case: dimension n for the index set

In most cases, for applications in other fields, results in dimension larger than one are required. Generalizing results obtained in dimension one seems often possible in theory, but, depending upon the method, it can reveal itself more difficult than expected.

There are three ways to consider the multidimensional case; in the first, the stochastic process is indexed by a set of dimension larger than one, in the second the process is such that each of its components is a random vector, and in the third these two previous cases are combined.

We will be interested in the first case, on which Wschebor (see [163]) has been working and generalizing results of Cramér et al. and Marcus for a *d*-parameter stochastic process $X = (X_t, t \in I\!\!R^d)$ having $C^1$ paths.

Wschebor studied problems related to the level sets of the paths of $X$ defined by $C_u^X := \{t : X_t = u\}$, proving at the same time some type of Rice formula. For more detail, see his 1985 book.

It is in this context that the work of Kratz and León (see [82]) on Gaussian fields took place. They considered the problem of the asymptotic behavior of the length of a level curve of a Gaussian field, by adapting the method used in the one dimensional case, for random processes indexed by a set of dimension larger than one.

Note that once again this study was motivated because of its various applications, in particular for the random modelisation of the sea (see e.g. [9]).

Consider a mean zero stationary Gaussian randon field $(X_{s,t}; (s,t) \in I\!\!R^2)$ with variance one and correlation function $r$ having partial derivatives $\partial_{ij} r$,
for $1 \le i + j \le 2$.

Assume that $r \in L^1$ satisfies $r^2(0,0) - r^2(s,t) \neq 0$, for $(s,t) \neq (0,0)$, and that $\partial_{02} r$, $\partial_{20} r$ are both in $L^2$.

Let $H(X)$ be the space of real square integrable functionals of the field $(X_{s,t}; (s,t) \in I\!\!R^2)$.

Let $\mathcal{L}_{Q(T)}^X(u)$ be the length of $\{(s,t) \in Q(T) : X_{s,t} = u\}$ the level curve at level $u$ for the random field $X$, $Q(T)$ being the square $[-T,T] \times [-T,T]$ and $|Q(T)|$ its Lebesgue measure.

By theorem 3.2.5 of Federer (see [53], p.244), we have for $g \in \mathbf{C}(I\!\!R)$

$$\int_{-\infty}^{\infty} g(u) \mathcal{L}_{Q(T)}^X(u) du = \iint_{Q(T)} g(X_{s,t}) ||\nabla X_{s,t}|| ds dt.$$

We assume w.l.o.g. that $\lambda := I\!\!E[\partial_{10}^2 X_{s,t}] = I\!\!E[\partial_{01}^2 X_{s,t}] = 1$ and that $X$ is isotropic, i.e. that the covariance matrix of $(\partial_{10} X_{s,t}, \partial_{01} X_{s,t})$ is of the form $\begin{pmatrix} \lambda & 0 \\ 0 & \lambda \end{pmatrix}$ (see [1], §6.2).

(Note that in the non-isotropic case, i.e. when the density function of $(\partial_{10} X_{s,t}, \partial_{01} X_{s,t})$ is $\mathcal{N}(0, \Sigma)$ with $\Sigma = \begin{pmatrix} \Sigma_{11} & \Sigma_{12} \\ \Sigma_{12} & \Sigma_{22} \end{pmatrix}$, $\Sigma_{i,j} \neq 0$ for $i \neq j$, we can consider the isotropic process $Y_{s,t}$ defined by $Y_{s,t} = \frac{1}{\sqrt{r(0,0)}} X_{u,v}$, where $\begin{pmatrix} u \\ v \end{pmatrix} = \frac{1}{\sqrt{r(0,0)}} \Sigma^{-1/2} \begin{pmatrix} s \\ t \end{pmatrix}$, and deduce the results for $X_{s,t}$ from the ones of $Y_{s,t}$.)

Under these conditions, we can obtain the chaos expansion in $H(X)$ for $\mathcal{L}_{Q(T)}^X(u)$, as well as a CLT for it, namely:

**Theorem 3.3** *(Kratz and León, 2001)*
*In $H(X)$, we have*

$$\mathcal{L}_{Q(T)}^X(u) = \sum_{q=0}^{\infty} \sum_{0 \le m+l \le [q/2]} d_{q-2(m+l)}(u) c_{2m,2l} \times$$

$$\iint_{Q(T)} H_{q-2(m+l)}(X_{s,t}) H_{2m}(\partial_{10} X_{s,t}) H_{2l}(\partial_{01} X_{s,t}) ds dt$$

*with $d_k(u) = \frac{1}{k!} H_k(u)\phi(u)$ ($\phi$ being the standard normal density) and*

$c_{2m,2l} = \dfrac{(-1)^{m+l}\sqrt{2\pi}}{m!l!2^{m+l}} \displaystyle\sum_{p_1=0}^{l} \sum_{p_2=0}^{m} \binom{l}{p_1}\binom{m}{p_2} \dfrac{(-1)^{p_1+p_2}}{B(p_1+p_2+1,1/2)}$ *(B being the Beta function).*

*The asymptotical behavior of $\mathcal{L}_{Q(T)}^X(u)$ is described by*

$$\frac{\mathcal{L}_{Q(T)}^X(u) - I\!\!E[\mathcal{L}_{Q(T)}^X(u)]}{|Q(T)|^{1/2}} \xrightarrow[T\to\infty]{} \mathcal{N}(0,\sigma^2),$$

*with $I\!\!E[\mathcal{L}_{Q(T)}^X(u)] = |Q(T)|\dfrac{\sqrt{2\pi}\phi(u)}{B(1,1/2)}$ and*

$$\sigma^2 = \sum_{q=0}^{\infty} \sum_{0 \le m_1+l_1 \le [q/2]} \sum_{0 \le m_2+l_2 \le [q/2]} \delta_{q,2m_1,2l_1}(u)\delta_{q,2m_2,2l_2}(u) \times$$

$$\iint_{I\!\!R^2} I\!\!E[I_{q,2m_1,2l_1}(0,0) I_{q,2m_2,2l_2}(s,t)] ds dt,$$

*where $\delta_{q,2m,2l}(u) := d_{q-2(m+l)}(u) c_{2m,2l}$, and*
*$I_{q,2m,2l}(s,t) := H_{q-2(m+l)}(X_{s,t}) H_{2m}(\partial_{10} X_{s,t}) H_{2l}(\partial_{01} X_{s,t})$.*

The proof is an adaptation to dimension 2 of the methods used to obtain proposition 2.2 and theorem 2.26 respectively.

Indeed we introduce $\mathcal{L}_{Q(T)}^X(u,\sigma) := \dfrac{1}{\sigma}\displaystyle\int_{I\!\!R} \mathcal{L}_{Q(T)}^X(v)\phi\left(\dfrac{u-v}{\sigma}\right) dv$ and prove that it converges to $\mathcal{L}_{Q(T)}^X(u)$ in $L^2$. Then the generalization to dimension 2 of lemma 2.6 will give us the Hermite expansion of $\mathcal{L}_{Q(T)}^X(u,\sigma)$ with coefficients $d_k^\sigma(u) c_{2m,2l}$, with $d_k^\sigma(u) \xrightarrow[\sigma\to 0]{} d_k(u)$.

**Lemma 3.1** *(Kratz and León, 2001)*
*Let $f \in L^2(\phi(u,v)dudv)$ and let $(d_k, k \ge 0)$ be its Hermite coefficients.*
*One has the following expansion*

$$\iint_{Q(T)} f(X_{s,t})||\nabla X_{s,t}||ds dt =$$

$$\sum_{k=0}^{\infty} \sum_{m=0}^{\infty} \sum_{l=0}^{\infty} d_k c_{2m,2l} \iint_{Q(T)} H_k(X_{s,t}) H_{2m}(\partial_{10} X_{s,t}) H_{2l}(\partial_{01} X_{s,t}) ds dt =$$

$$\sum_{k=0}^{\infty} \sum_{0 \leq m+l \leq [q/2]} d_{k-2(m+l)} c_{2m,2l} \int\int_{Q(T)} H_{k-2(m+l)}(X_{s,t}) H_{2m}(\partial_{10} X_{s,t}) H_{2l}(\partial_{01} X_{s,t}) ds dt,$$

*where $c_{2m,2l}$ is given in theorem 3.3.*

The first part of theorem 3.3, i.e. the chaos expansion follows then, exactly as in dimension one.

As regards the CLT, we define

$$L_{Q(T)}(u) = \frac{\mathcal{L}_{Q(T)}^X(u) - I\!E[\mathcal{L}_{Q(T)}^X(u)]}{|Q(T)|^{1/2}}$$

and $L_{Q(T)}^K$ as the finite sum deduced from $L_{Q(T)}(u)$ for $q = 1$ to $K$, i.e.

$$L_{Q(T)}^K(u) = \frac{1}{|Q(T)|^{1/2}} \sum_{q=1}^{K} \sum_{0 \leq m+l \leq [q/2]} \delta_{q,2m,2l}(u) \int\int_{Q(T)} I_{q,2m,2l}(s,t) ds dt.$$

We define in the same way $L_{Q(T)}^K(u, \varepsilon)$ in which we consider the $1/\varepsilon$-dependent random field $X^\varepsilon$ defined as

$$X_{s,t}^\varepsilon = \int_{I\!R^2} e^{i(s\lambda_1 + t\lambda_2)} ((f * \hat{\varphi}_\varepsilon)(\lambda_1, \lambda_2))^{\frac{1}{2}} dW(\lambda_1, \lambda_2).$$

We must however normalize the partial derivatives, dividing by the constant $(I\!E[\partial_{10} X_{0,0}^\varepsilon]^2)^{1/2} = (I\!E[\partial_{01} X_{0,0}^\varepsilon]^2)^{1/2} = (-\partial_{20} r(0,0))^{1/2}$ to get random variables with variance one.

The proof follows then the one of theorem 2.26, using our combined method.

With respect to Rice formulas and the distribution of the maximum of Gaussian fields, let us also mention the recent works by Adler et al. and Taylor et al. (see [157] and [2]), and by Azaïs and Wschebor (see [14]).

## 4. Conclusion

Although this work on level crossings focuses specifically on the stationary Gaussian case, it is to be noted that much research has also been conducted on the non-stationary Gaussian case (notion of local stationarity, notion of non-stationarity but with constant variance; work on diffusions, Brownian motions, fractional Brownian motions, etc...), as well as on the non-Gaussian case (in particular when considering the class of stable processes).

For completeness, we tried to make a large inventory of papers and books dealing with the subject of level crossings, even though all the references are not explicitly mentioned in the synopsis above.

Note that in case an author published a paper prior to a book, the only reference mentioned is that of his/her book.

Marie Kratz
SID Department
ESSEC
av. Bernard Hirsch, BP50105
95021 Cergy-Pontoise cedex
FRANCE
`http://www.essec.fr/professeurs/marie-kratz`
or `http://www.math-info.univ-paris5.fr/~kratz`